\documentclass{amsart}
\usepackage{amssymb}
\usepackage{euscript}
\usepackage{amsfonts}
\usepackage{epsfig}
\usepackage{color}
\usepackage[frame,ps,dvips,matrix,arrow,curve,rotate]{xy}
\let\oldcite\cite
\newtheorem{thm}{Theorem}[section]
\newtheorem{cor}[thm]{Corollary}
\newtheorem{lem}[thm]{Lemma}

\newtheorem{prop}[thm]{Proposition}
\theoremstyle{definition}
\newtheorem{defn}[thm]{Definition}
\theoremstyle{remark}
\newtheorem{rem}[thm]{Remark}
\numberwithin{equation}{section}
\theoremstyle{remark}
\newtheorem{ex}[thm]{Example}

\newcommand{\bbG}{\mathbb{G}}
\newcommand{\bbH}{\mathbb{H}}
\newcommand{\bbF}{\mathbb{F}}
\newcommand{\bbK}{\mathbb{K}}
\newcommand{\mfC}{\mathfrak{C}}
\newcommand{\mfD}{\mathfrak{D}}
\newcommand{\mfA}{\mathfrak{A}}
\newcommand{\M}{\mathcal{B}}

\newcommand{\A}{\mathsf{A}}
\newcommand{\Ch}{\mathcal{C}h^{[-1,0]}}
\newcommand{\Chst}{\mathcal{C}h^{[-1,0]}_{st}}

\newcommand{\cm}{\mathbf{XMod}}
\newcommand{\tgp}{\mathbf{2Gp}}
\newcommand{\tgpd}{\mathbf{2Gpd}}

\let\:=\colon

\newcommand{\lra}{\longrightarrow}
\newcommand{\llra}[1]{\stackrel{#1}{\lra}}

\newcommand{\Ob}{\operatorname{Ob}}

\newcommand{\Out}{\operatorname{Out}}
\newcommand{\Aut}{\operatorname{Aut}}

\newcommand{\Hom}{\operatorname{Hom}}
\newcommand{\homo}{\underline{\mathcal{H}om}_*}
\newcommand{\homc}{\underline{\mathcal{H}om}}
\newcommand{\Hombf}{\mathbf{Hom}}
\newcommand{\rhom}{\underline{\mathcal{RH}om}}
\newcommand{\Map}{\underline{\mathcal{RH}om}_*}
\newcommand{\hhom}{\Hom_{\Ho(\mathbf{2Gp})}}
\newcommand{\homw}{\underline{\mathcal{HOM}}}

\newcommand{\id}{\operatorname{id}}
\newcommand{\Ho}{\operatorname{Ho}}
\newcommand{\Ker}{\operatorname{Ker}}
\newcommand{\Coker}{\operatorname{Coker}}

\def\smashedlongrightarrow{\setbox0=\hbox{$\longrightarrow$}\ht0=1pt\box0}
\def\risom{\buildrel\sim\over{\smashedlongrightarrow}}
\def\smashedst{\setbox0=\hbox{$\rightrightarrows$}\ht0=4pt\box0}
\def\smashedar{\setbox0=\hbox{$\longrightarrow$}\ht0=4pt\box0}
\newcommand{\sst}[1]{\stackrel{#1}{\smashedst}}
\newcommand{\arr}[1]{\stackrel{#1}{\smashedar}}
\newcommand{\oux}[2]{\underset{#1}{\overset{#2}\times}}

\begin{document}

\title[On weak maps between 2-groups]
      {On weak maps between 2-groups}%

\author{Behrang Noohi}%

\begin{abstract}
  We give an explicit handy  cocycle-free description of the
  groupoid of weak maps between two crossed-modules using what we
  call a {\em butterfly} (Theorem \ref{T:butterfly}). We define
  composition of butterflies and this way find a bicategory that is
  naturally biequivalent to the 2-category of pointed homotopy
  2-types. This has applications in the study of 2-group actions
  (say, on stacks), and in the theory
  of gerbes bound by crossed-modules and principal-2-bundles).
\end{abstract}
\maketitle
\section{Introduction}{\label{S:Intro}}

There are several incarnations of 2-groups in mathematics. To
mention a few:\footnote{To see some more explicit example see the
last section of \oldcite{Baez}.} 1) They appear abstractly as
special classes of monoidal categories, e.g., in the context of
bitorsors \cite{Breen} or Picard categories \cite{Del}; 2) They
classify connected pointed homotopy 2-types (via the fundamental
2-group), as discovered by MacLane and Whitehead; 3) They appear as
symmetries of objects in a 2-category. Let us dissect  (3) a bit
further by giving some examples.

The auto-equivalences of a stack from a 2-group. Thus, 2-groups play
an important rule in the study of (2-) group actions on stacks
\cite{B-N}. Along similar lines, 2-groups appear in the
representation theory of 2-vector spaces \cite{El}, \cite{GaKa}, and
also in the theory of gerbes  and higher local systems \cite{Al},
\cite{Deb}, [Bre1-4], \cite{BrMe}, \cite{PoWa}, \cite{BDR}. The
latter is indeed deeply connected with physics through higher gauge
theory: various kind of gerbes that come up in gauge theory are to
be thought of as 2-principal-bundles for a certain 2-group. For
example: Giraud's $G$-gerbes are 2-principal-bundles for the 2-group
$\mathfrak{Aut}(G)$ \cite{Breen}; $S^1$-gerbes \cite{Bry} are
2-principal-bundles for the 2-group whose 2-morphisms are $S^1$ and
has no nontrivial 1-morphisms; string bundles are essentially the
same as 2-principal-bundles for the string 2-groups \cite{BCSS},
\cite{BaSch}. Loop groups are also closely related to 2-groups
\cite{BCSS}.

Most 2-groups that arise in nature are {\em weak}, in the sense
that, either the multiplication is only associative up to higher
coherences, or inverses exist only in a weak sense (or both).
Furthermore, interesting morphisms between 2-groups are also often
weak, in the sense that they preserve products only up to higher
coherences. It is a standard fact that one can strictify weak
2-groups, but not the weak functors. Put differently, the 2-category
of {\em strict} 2-groups and {\em weak} functors between them is the
``homotopically correct'' habitat for 2-groups.

The (strict) 2-groups by themselves can be codified conveniently
using crossed-modules. Weak morphisms between strict 2-groups,
however, are more complicated and to write them down results in
somewhat inconvenient cocycles. It is therefore desirable to find a
clean and way to deal with weak morphisms between 2-groups so as to
make them more tractable in geometric situations.

The aim of this paper is to do exactly this. Namely, we give a
concrete and manageable cocycle-free model for the space of weak
morphisms between two crossed-modules. (It is easy to see that this
space is a 1-type, so its homotopy type is described by a groupoid.)
{\bf Theorem \ref{T:butterfly}} (also see {\bf Theorem
\ref{T:classification}}) gives us a functorial model for this
groupoid in terms of what we call {\em butterflies}. Butterflies
indeed furnish a neat bicategory structure on crossed-modules ({\bf
Theorem \ref{T:bicat}}), therefore giving rise to a model for the
homotopy category of pointed connected 2-types. We also discuss the
braided and abelian versions of this bicategory.

In a future paper we  generalize these result to the case where
everything is relative to a Grothendieck site. As a consequence, the
main results of the present paper will apply to the case of Lie
2-groups, 2-group schemes, and so on. Bear in mind that in these
geometric settings the cocycle approach to weak morphisms is very
inconvenient and sometimes hopeless (see below).

\vspace{01in}\noindent{\bf Some applications.}

For an application of butterflies to the HRS-tilting theory
\cite{HaReSm} we refer the reader to \cite{Tilting}. In another
application (joint work with E.~Aldrovandi), we investigate
butterflies from the point of view of gr-stacks. We also employ
butterflies (over a Grothendieck site) to study the ``change of the
structure 2-group along a weak morphism of 2-groups'' for
2-principal-bundles and its effect on their geometry.
 (This has been previously looked at only for strict 2-group
morphism \cite{Breen}.)

\vspace{0.1in} \noindent{\bf 2-group actions on stacks.} Let us
spell out in more detail an application of our approach to weak
morphisms -- this will also serve to explain why a cocycle-free
approach could be advantageous sometimes.

We propose a systematic way to study 2-group actions on stacks.
Given a stack $\mathcal{X}$ (say, topological, differentiable,
analytic, algebraic, etc.), the set
 $\Aut \mathcal{X}$ of self-equivalences of $\mathcal{X}$
is naturally  a weak 2-group. (It is weak because equivalences are
not strictly invertible; the associativity, however, remains
strict.) To have a 2-group $\bbH$ act on $\mathcal{X}$ is the same
thing as to have a weak map from $\bbH$ to $\Aut \mathcal{X}$. If
two such maps are related by a (pointed) transformation, they should
be regarded as giving the ``same'' action of $\bbH$ on
$\mathcal{X}$. So the question is to classify such equivalence
classes of weak maps $f \: \bbH \to \Aut \mathcal{X}$.

The weakness of $\Aut \mathcal{X}$, and of the map $f$, are,
however, disturbing and we would like to make things as strict as
possible. Using a bit of homotopy theory, and some standard
strictification procedures, it can be shown that what we are looking
for is $[\bbH,\bbG]_{\tgp}$, where $\bbG$ is a strict model for
$\Aut \mathcal{X}$. Here, $[\bbH,\bbG]_{\tgp}$ stands for the set of
morphisms from $\bbH$ to $\bbG$ in the homotopy category $\Ho(\tgp)$
of the category of strict 2-groups and strict maps. Equivalence of
2-groups and crossed-modules ($\S$ \ref{SS:Equiv}) now enables us to
translate the problem to the language of crossed-modules, in which
case we have an explicit description of $[\bbH,\bbG]_{\tgp}$ in
simple group theoretic terms thanks to Theorem \ref{T:butterfly}
(also see Theorem \ref{T:classification} and Corollary
\ref{C:homotopy}).

All we need is to run this method is to find a crossed-module  model
for $\Aut \mathcal{X}$

Thanks to the very explicit nature of the above procedure, we are
 able to give solid constructions with stacks, circumventing a
lot of ``weaknesses'' and coherence conditions that arise in
studying group actions on stacks. An application of this
strictification method is given in \cite{B-N}, where the covering
theory of stacks is used to classify smooth Deligne-Mumford analytic
curves, and also to give an explicit description of them as quotient
stacks. For instance, using these 2-group theoretic techniques we
obtain the following completely geometric result:  every smooth
analytic (respectively, algebraic) Deligne-Mumford stack of
dimension one is the quotient stack for the action of either a
finite group or a central finite extension of $\mathbb{C}^*$ on a
complex manifold (respectively, complex variety).

 \vspace{0.1in} \noindent{\bf Organization  of the paper}

\vspace{0.1in}

Sections \ref{S:Review} to \ref{S:M-S} are devoted to recalling some
standard facts about 2-groups and crossed-modules and fixing the
notation. Essential for more easily reading the paper is the fact
that the category of 2-groups is equivalent to the category of
crossed-modules ($\S$\ref{SS:Equiv}). The reader will find it
beneficial to keep in mind how this equivalence works, as we will
freely switch back and forth between 2-groups and crossed-modules
throughout the paper (sometimes even using the two terms
synonymously). For us, 2-groups are the {\em conceptual} side of the
story, whereas crossed-modules provide the {\em computational}
framework.

Viewed as 2-groupoids with one object, 2-groups (hence, also
crossed-modules) can be treated via the Moerdijk-Svensson model
structure \cite{M-S}. This is briefly recalled in Section
\ref{S:M-S}. We point out that, all we need from closed model
categories is the notion of fibrant/cofibrant resolution and the way
it can be used to compute hom-sets in the homotopy category. Taking
this for granted, the reader unfamiliar with closed model categories
can proceed with no difficulty.

Section \ref{S:group} concerns some elementary  constructions from
group theory.  We introduce a pushout construction for
crossed-modules and work out its basic properties. This section is
perhaps is not so interesting by itself, but it provides the
technical tools required in the proof of Theorem \ref{T:butterfly}.

Section \ref{S:Butterfly} is the core of the paper. In it we state
and prove our main result (Theorem \ref{T:butterfly}). This is based
on the notion of {\em butterfly} (Definition \ref{D:butterfly}). We
investigate butterflies some more in Section \ref{S:More}. In
$\S$\ref{S:Bicat} we give an explicit model $\mathcal{CM}$ for the
2-category of crossed-modules and weak morphisms in terms of
butterflies.

In Section \ref{S:Ker} we indicate how the notions of kernel,
cokernel, exact sequence, and so on of weak morphisms of 2-groups
find a natural and simple form in the world of butterflies.

Braided and  abelian butterflies are briefly discussed in Section
\ref{S:Derived}. We use the latter to give a simple  description of
the derived category of complexes of length 2 in an abelian category
$\A$.

In Section \ref{S:special} we consider a special case of Theorem
\ref{T:butterfly} in which the source 2-group is an honest group
(Theorem \ref{T:classification}) and discuss its connection with
results of Dedecker and Blanco-Bullejos-Faro. In Section
\ref{S:Cohomological} we discuss a cohomological version of this
(Theorem \ref{T:maps}). This cohomological classification is not
 a new result (with some diligence, the reader can verify that it is
a special case of the work of \cite{AzCe}), and our emphasis is only
to make precise the way it relates to Theorem
\ref{T:classification}, as it was used in \cite{B-N}.

In Sections \ref{S:Homotopical}  and \ref{S:Compatible} we explain
the homotopical meaning of Theorem \ref{T:classification} in terms
the Postnikov decomposition of the classifying space of a
crossed-module; again, this is folklore. We make this precise using
the notion of {\em difference fibration} which is an obstruction
theoretic construction introduced in \cite{Baues} used in studying
the liftings of  a map  into a fibration (from the base to the total
space).

In the appendix we review basic general facts about 2-categories
and 2-groupoids.

\vspace{0.1in} \noindent{\em Acknowledgement.}  I am grateful to
E.~Aldrovandi, H-J.~Baues, L.~Breen, J.~Elgueta, M.~Jibladze,
M.~Kamgarpour, and F.~Muro for many useful discussions. I would like
to thank Bertrand To\"en for pointing me to Elgueta's work
\cite{El}. The idea of Theorem \ref{T:classification} was conceived
in conversations with Kai Behrend in our joint work \cite{B-N}. The
relation between butterflies and cocycles for weak morphisms was
pointed out to me by Roman Mikhailov. Finally, I would like to thank
Max-Planck-Institut f\"{u}r Mathematik, Bonn, where the research for
this work was done, for providing pleasant working conditions.

\tableofcontents

\section{Notation and terminology}{\label{S:NT}}

We list some of the notations and conventions used throughout the
paper.

The structure map of a crossed-module $\bbG=[G_2\to G_1]$ is usually
denoted by $\alpha \mapsto \underline{\alpha}$. The components of a
morphism $P: \bbH \to \bbG$ of crossed-modules are denoted by $p_2
\: H_2 \to G_2$ and $p_1 \: H_1 \to G_1$. The  action of $G_1$ on
$G_2$ is denoted by $-^a$, and so is the conjugation action of $G_1$
on itself.

In a semi-direct product $A\rtimes B$, respectively $B\ltimes A$,
the group $B$ acts on $A$ on the left, respectively right.

For objects $A$ and $B$ in a  category $\mathbf{C}$ with a notion of
weak equivalence, we denote the the set  of morphisms in the
homotopy category from $B$ to $A$, that is
$\Hom_{\Ho(\mathbf{C})}(B,A)$, by $[B,A]_{\mathbf{C}}$.

We usually  denote a short exact sequence
    $$ 1 \to N \to E \to \Gamma \to 1$$
simply by $E$. Quotient maps, such as $E \to \Gamma$ in the above
sequence, or $G_1 \to \pi_1\bbG$, are usually denoted by $x \mapsto
\bar{x}$.

\vspace{0.1in} \noindent{\bf Some abuse of terminology.}

We tend to use the term ``map'' where it is perhaps more appropriate
to use the term ``morphism''. Also, we use the term ``equivalence''
(of 2-groups, crossed-modules, etc.) for what should really be
called a ``weak equivalence'' (or ``quasi-isomorphism'').

For 2-groups $\bbG$ and $\bbH$, when we say  {\em the homotopy class
of a weak map from $\bbH$ to $\bbG$}, we mean a map in $\Ho(\tgp)$
from $\bbH$ to $\bbG$; this terminology is justified by Theorem
\ref{T:weak} and Proposition \ref{P:mapping}.

\section{Quick review of 2-groups and crossed-modules}{\label{S:Review}}

\subsection{Quick review of 2-groups}{\label{SS:2gp}}

We recall some basic facts about 2-groups and crossed-modules. Our
main references are \cite{Brown,W,M-W,M-S,Loday,Baez}.

A {\bf 2-group} $\mathfrak{G}$ is a  group object in the category of
groupoids. Alternatively, we can define a 2-group to be a groupoid
object in the category of groups, or  also, as a (strict) 2-category
with one object in which all  1-morphisms and 2-morphisms are
invertible (in the strict sense). We will try to stick with the
`group in groupoids' point of view throughout the paper, but
occasionally switching back and forth between different points of
view is inevitable. Therefore, the reader will  find it rewarding to
master how the equivalence of these three point of views works.

A (strict) {\em morphism} $f \: \mathfrak{G} \to \mathfrak{H}$ of
2-groups is a map of groupoids that respects the group operation. If
we view $\mathfrak{G}$ and  $\mathfrak{H}$ as 2-categories with one
object, such $f$ is   nothing but  a strict 2-functor. The category
of 2-groups is denoted by $\tgp$.

To a 2-group $\mathfrak{G}$ we associate the groups
$\pi_1\mathfrak{G}$ and $\pi_2\mathfrak{G}$ as follows. The group
$\pi_1\mathfrak{G}$ is the set of isomorphism classes of object of
the groupoid $\mathfrak{G}$. The group structure on
$\pi_1\mathfrak{G}$ is induced from the group operation of
$\mathfrak{G}$. The group $\pi_2\mathfrak{G}$ is the group of
automorphisms  of the identity object $e \in \mathfrak{G}$. This is
an abelian group. A  morphism $f \: \mathfrak{G} \to \mathfrak{H}$
of 2-groups is called an {\em equivalence} if it induces
isomorphisms on $\pi_1$ and $\pi_2$. The {\em homotopy category} of
2-groups is the category obtained by inverting all the equivalences
in $\tgp$. We denote it by $\Ho(\tgp)$.

\vspace{0.1in} \noindent{\em Caveat}: an equivalence between
2-groups need not have an inverse. Also, two equivalent 2-groups may
not be related by an equivalence, but only a zig-zag of
equivalences.

\subsection{Quick review of Crossed-modules}{\label{SS:CrossedMod}}

A {\bf crossed-module} $\mathfrak{G}=[G_2 \to G_1]$ is a pair of
groups $G_1,G_2$, a group homomorphism $\partial_{\bbG} \: G_2 \to
G_1$, and a (right) action of $G_1$ on $G_2$, denoted $-^a$, which
lifts to $G_2$ the conjugation action of $G_1$ on the image of
$\partial_{\bbG}$ and descends to $G_1$ the conjugation action of
$G_2$ on itself. The kernel of $\partial$ is a central (in
particular abelian) subgroup of $G_2$ and is denoted by
$\pi_2\mathfrak{G}$. The image of $\partial$ is a normal subgroup of
$G_1$ whose cokernel is denoted by $\pi_1\mathfrak{G}$. A  (strict)
morphism  of crossed-modules is a pair of group homomorphisms which
commute with the $\partial$ maps and respect the actions. A morphism
is called an {\em equivalence} if it induces isomorphisms on $\pi_1$
and $\pi_2$.

\vspace{0.1in} \noindent{\bf Notation.}   Elements of $G_2$ are
usually denoted by Greek letters and those of $G_1$ by lower case
Roman letters. The components of a map $P: \bbH \to \bbG$ of
crossed-modules are denoted by $p_2 \: H_2 \to G_2$ and $p_1 \: H_1
\to G_1$. We usually use $\partial$ instead of $\partial_{\bbG}$ if
the 2-group $\bbG$ is clear from the context. We sometimes suppress
$\partial$ from the notation and denote $\partial(\alpha)$ by
$\underline{\alpha}$.  For elements $g$ and $a$ in a group $G$ we
sometimes denote  $a^{-1}ga$ by $g^a$. The compatibility assumptions
built in the definition of a crossed-module make this unambiguous.
With this notation, the two compatibility axioms of a
crossed-modules can be written in the following way:
 \begin{itemize}
   \item[$\mathbf{CM1.}$] $\forall \alpha,\beta \in G_2, \
                                \beta^{\underline{\alpha}}=\beta^{\alpha}$;

   \item[$\mathbf{CM2.}$] $\forall \beta \in G_2, \forall a\in G_1, \
                                   \underline{\beta}^a=\underline{\beta^a}$.
 \end{itemize}

\subsection{Equivalence of 2-groups and crossed-modules}{\label{SS:Equiv}}

There is a natural pair of inverse equivalences   between the
category $\tgp$ of 2-groups  and the category $\cm$ of
crossed-modules. Furthermore, these functors preserve $\pi_1$ and
$\pi_2$. They are constructed as follows.

\vspace{0.1in}
\noindent{\em Functor from 2-groups to
crossed-modules.} Let $\mathfrak{G}$ be a 2-group. Let $G_1$ be the
group of   objects of $\mathfrak{G}$, and $G_2$  the set of   arrows
emanating from the identity object $e$; the latter is also a group
(namely, it is a subgroup of the group  of arrows of
$\mathfrak{G}$).

Define the map $\partial \: G_2 \to G_1$ by sending $\alpha \in G_2$
to $t(\alpha)$.

The action of $G_1$ on $G_2$ is given by conjugation. That is, given
$\alpha \in G_2$ and $g \in G_1$, the action is given by
$g^{-1}\alpha g$. Here were are thinking of $g$ as an identity arrow
and multiplication takes place in the group of arrows of of $\bbG$.
It is readily checked that $[\partial \: G_2 \to G_1]$ is a
crossed-module.

\vspace{0.1in} \noindent{\em Functor from crossed-modules to
2-groups.} Let $[\partial \: G_2 \to G_1]$ be a crossed-module.
Consider the groupoid $\mathfrak{G}$ whose underlying set of objects
is $G_1$ and whose set of arrows is $G_1\ltimes G_2$. The source and
target maps are given by $s(g,\alpha)=g$,
$t(g,\alpha)=g\partial(\alpha)$. Two arrows $(g,\alpha)$ and
$(h,\beta)$ such that $g\partial(\alpha)=h$ are composed to
$(g,\alpha\beta)$. Now, taking into account the group structure on
$G_1$ and the semi-direct product group structure on $G_1\ltimes
G_2$, we see that $\mathfrak{G}$ is indeed  a groupoid  object in
the category of groups, hence a 2-group.

\vspace{0.1in} The above discussion shows that there is a pair of
inverse functors inducing an equivalence between $\cm$ and $\tgp$.
These functors respect $\pi_1$ and $\pi_2$. Therefore, we have an
equivalence
   $$\xymatrix@=16pt@M=8pt{
      \Ho(\cm)  \ar@<0.5ex>[r]
                 & \ar@<0.5ex>[l] \Ho(\tgp).}$$

\section{Transformations between morphisms of crossed-modules}{\label{S:Transfo}}

We go over the notions of transformation and pointed transformation
between maps of crossed-modules. These are adaptations of the usual
2-categorical notions, translated to the crossed-module language via
the equivalence $\cm\cong\tgp$ (see Appendix,
$\S$\ref{SA:Transformations}). The idea is to think of a
crossed-module as a 2-group ($\S$\ref{SS:Equiv}), which is itself
thought of as a 2-groupoid with one object.

 \begin{defn}{\label{D:transformation}}
  Let $\bbG=[G_2 \to G_1]$ and $\bbH=[H_2 \to
  H_1]$ be crossed-modules, and let $P,Q \: \bbH \to \bbG$ be
  morphisms between them.
  A {\em transformation}  $T \: Q \Rightarrow P$ consists of
  a pair $(a,\theta)$ where $a \in G_1$ and $\theta \: H_1 \to G_2$
  is a crossed homomorphism  for the induced action, via $p_1$, of
  $H_1$ on $G_2$ (that is, $\theta(hh')=\theta(h)^{p_1(h')}\theta(h')$).
  We require the following:

    \begin{itemize}
        \item[$\mathbf{T1.}$] $q_1(h)^a=p_1(h)\underline{\theta(h)}$, for
             every $h \in H_1$;
        \vspace{1mm}
        \item[$\mathbf{T2.}$] $q_2(\beta)^a=p_2(\beta)\theta(\underline{\beta})$,
             for every  $\beta \in H_2$.
    \end{itemize}
  We say $T$ is {\em pointed} if $a=1$; in this case, we denote $T$
  simply by $\theta$.    When $\theta$ is the trivial map,
  the transformation
  $T$ is called  {\em conjugation by} $a$; in this case, we use the notation
  $P=Q^a$ or $P=a^{-1}Qa$.
 \end{defn}

 \begin{rem}{\label{R:conj}}
  Given a 2-group $\mathfrak{G}$ and an element $a$ in $G_1$ (the
  group of objects) we define the  morphism $c_a \: \bbG \to \bbG$,
  called {\em conjugation by $a$},   to be the map that sends an
  object $g$ (respectively, an arrow $\alpha$)   to $a^{-1}ga$
  (respectively, $a^{-1}\alpha a$). If we consider the corresponding
  crossed-module $[G_2 \to G_1]$, the conjugation morphism
  $c_a$ sends  $g \in G_1$ to $a^{-1}ga$ and $\alpha \in G_2$ to
  $\alpha^a$. In the notation of Definition
  \ref{D:transformation}, it is easy to see that $Q^a=c_a\circ Q$.
 \end{rem}

 \begin{lem}{\label{L:conj}}
   Let $P,Q \: \bbH \to \bbG$ be maps of 2-groups and $(a,\theta)$ a
   transformation from $Q$ to $P$. Then $\pi_iP=(\pi_iQ)^a$,
   $i=1,2$. In particular, if $P$ and $Q$ are related by a pointed
   transformation, then they induce the
   same map on  homotopy groups.
 \end{lem}

  \begin{proof}
    Obvious.
  \end{proof}

Let $P,Q,R \: \bbH \to \bbG$ be maps of crossed-modules. Given
homotopies $(b,\sigma) \: R \Rightarrow Q$ and $(a,\theta) \: Q
\Rightarrow P$,  consider the pointwise product $\theta\sigma \: H_1
\to G_2$. It is easily checked that $(ba,\theta\sigma)$ is a
transformation from $R$ to $P$. 
This construction, of course, corresponds
to the usual composition of weak 2-transformation between
2-functors.

A transformation $(a,\theta) \: Q \Rightarrow P$ has an inverse
$(a^{-1},\theta^{-1}) \: P \Rightarrow Q$, where $\theta^{-1} \: H_1
\to G_2$ is defined by  $\theta^{-1}(h):=\theta(h)^{-1}$.

 \begin{defn}{\label{D:homo}}
   Let $\bbG$ and $\bbH$ be crossed-modules. We define the
   {\em mapping groupoid} $\homo(\bbH,\bbG)$
   to be the groupoid whose  objects are crossed-module maps
   $\bbH \to \bbG$ and whose morphisms are pointed transformations.
 \end{defn}

 \begin{rem} Observe that in the definition above of
     $\homo(\bbH,\bbG)$ we have not used `modifications' and,
     in particular, the outcome is a groupoid and not a
     2-groupoid. This is because between two pointed transformations
     there is no non-trivial {\em pointed} modification.
 \end{rem}

 With hom-groupoids being $\homo(\bbH,\bbG)$, the category $\cm$ is
 enriched over groupoids. We denote the resulting 2-category by $\underline{\cm}$.
 Similarly, $\tgp$ can be enriched over groupoids by taking
 2-morphisms to be pointed transformations (Definition
 \ref{D:transformation}). We denote the resulting 2-category by $\underline{\tgp}$.
 The equivalence
       $$\cm  \cong \tgp$$
 of $\S$\ref{SS:Equiv} now becomes a biequivalence of 2-categories

 \begin{prop}{\label{P:functoriality}}%
     The construction of $\S$\ref{SS:Equiv} gives rise to a
     biequivalence
       $$\underline{\cm}  \cong \underline{\tgp}.$$
     of 2-categories.
 \end{prop}

The biequivalence of Proposition \ref{P:functoriality} is a
biequivalence in a strong sense: it induces {\em isomorphisms} on
hom-groupoids.

 The mapping space $\homo(\bbH,\bbG)$ in not the homotopically ``correct''
mapping space, as it lacks the expected homotopy invariance
property. That is, an equivalence $\bbH' \to \bbH$ of
crossed-modules does not necessarily induce an equivalence
$\homo(\bbH,\bbG) \to \homo(\bbH',\bbG)$ of  groupoids. We explain
in $\S$\ref{S:M-S} how this failure can be fixed by making use of
cofibrant replacements in the category of crossed-modules
(especially, see Definition \ref{D:mapping}).

\section{Weak morphisms between 2-groups}{\label{S:Weak}}

A 2-group $\bbG$ can equivalently be defined to be  a strict4
monoidal groupoid in which multiplication (on the left, and on the
right) by any object is an isomorphism of categories. Given two
2-groups $\bbG$ and $\bbH$, the mapping groupoid $\homo(\bbH,\bbG)$
is then naturally isomorphic to the  groupoid whose object are
strict monoidal functors and whose morphisms are natural monoidal
transformations.

We define a {\bf weak morphism} $f \: \bbH \to \bbG$ to be a
monoidal functor which respects the unit objects strictly (also see
\cite{Notes}, $\S$7 and $\S$8). With monoidal transformations
between them, these are objects of a groupoid which we denote by
$\homw_*(\bbH,\bbG)$.\footnote{For the interpretation of  weak
morphisms in the crossed-module language see \oldcite{Notes},
$\S$8.}

The goal of the paper is to give an explicit model for
$\homw_*(\bbH,\bbG)$. Our strategy is to use the fact that
$\homw_*(\bbH,\bbG)$ is equivalent to the  derived mapping groupoid
$\Map(\bbH,\bbG)$; see Theorem \ref{T:weak}. We then give an
explicit model for the derived mapping space $\Map(\bbH,\bbG)$  in
terms of {\em butterflies} (Definition \ref{D:butterfly} and Theorem
\ref{T:butterfly}).

\section{Moerdijk-Svensson closed model structure
and crossed-modules}{\label{S:M-S}}

It has been known since \cite{W} that crossed-modules model pointed
connected homotopy 2-types. That is, the pointed homotopy type of a
connected pointed CW-complex with $\pi_iX=0$, $i \geq 3$, is
determined by (the equivalence class of) a crossed-module. In
particular, the homotopical invariants of such a CW-complex can be
read off from the corresponding crossed-module.

The approach in \cite{W} and \cite{M-W} to the classification of
2-types is, however, not functorial. To have a functorial
classification of homotopy 2-types (i.e. one that also accounts for
maps between such objects), it is best to incorporate  closed model
categories.  To do so, recall that a crossed-module can be regarded
as a 2-group, and a 2-group is in turn a 2-groupoid with one object.

In \cite{M-S}, Moerdijk and Svensson introduce a   closed model
structure  on the (strict) category of (strict) 2-groupoids, and
show that there is a Quillen pair between the closed model category
of 2-groupoids and the closed model category of CW-complexes, which
induce an equivalence between the homotopy category of 2-groupoids
and the homotopy category of CW-complexes with vanishing $\pi_i$,
$i\geq 3$. We use this model structure to deduce some results about
crossed-modules.

We emphasize that, in working with crossed-modules, what we are
using is the {\em pointed} homotopy category. So we need to adopt a
pointed version of the Moerdijk-Svensson structure. But this does
not cause any additional difficulty as everything in \cite{M-S}
carries over to the pointed case. For a quick review of the
Moerdijk-Svensson structure see Appendix.

It is easy to see that a weak equivalence between 2-groups in the
sense of Moerdijk-Svensson is the same as a weak equivalence between
crossed-modules in the sense of $\S$\ref{S:Review}.  Let us see what
the fibrations look like.

 \begin{defn}{\label{D:fibrations}}
   A map $(f_2,f_1) \: [H_2 \to H_1] \to  [G_2 \to G_1]$ of crossed-modules
   is called a {\em fibration} if $f_2$ and $f_1$ are both surjective.
   It is called {\em a trivial fibration} if, furthermore, the
   map $H_2 \to H_1\times_{G_1} G_2$ is an isomorphism.
 \end{defn}

We leave it to the reader to translate these to the language of
2-groups and verify that they coincide with Moerdijk-Svensson
definition of (trivial) fibration.

Let us now look at cofibrations. In fact, we will only describe what
the cofibrant objects are, because that is all we need in this
paper.

 \begin{defn}{\label{D:cofibrations}}
    A crossed-module $[G_2 \to G_1]$ is {\em cofibrant} if $G_1$
    is a free group.
 \end{defn}

Observe that this is much weaker than Whitehead's notion of a free
crossed-module. However, this is the one that corresponds to
Moerdijk-Svensson's definition.

 \begin{prop}
  A crossed-module $\mathfrak{G}=[G_2 \to G_1]$ is cofibrant in the
  sense of Definition \ref{D:cofibrations} if and only if its
  corresponding 2-group is cofibrant in the
  Moerdijk-Svensson structure.
 \end{prop}

 \begin{proof}
  This follows immediately from the Remark on page 194 of
  \cite{M-S}, but we give a direct proof. A 2-group $\bbG$ is
  cofibrant in Moerdijk-Svensson structure, if and only if every
  trivial fibration $\bbH \to \bbG$, where $\bbH$ is a 2-group{\em
  oid}, admits a section. But, we can obviously restrict ourselves
  to 2-groups $\bbH$. So, we can work entirely   within
  crossed-modules, and  use the notion of trivial fibration as
  in Definition \ref{D:fibrations}.

  Assume $G_1$ is free. Let $(f_2,f_1)\: [H_2 \to H_1] \to  [G_2 \to
  G_1]$ be a trivial fibration. Since $G_1$ is free and $f_1$ is
  surjective, there is a section $s_1 \: G_1 \to H_1$. Using  the
  fact that $H_2\cong H_1\times_{G_1} G_2$, we also get a natural
  section $s_2 \: G_2 \to H_2$ for the projection
  $H_1\times_{G_1}G_2\to G_2$, namely,
  $s_2(\alpha)=(s_1\big(\underline{\alpha}),\alpha\big)$. It is easy
  to see that $(s_2,s_1) \: [G_2 \to G_1] \to [H_2 \to H_1]$
  is a map of crossed-modules.

  To prove the converse, choose a free group $F_1$ and a surjection
  $f_1 \: F_1 \to G_1$. Form the pull back crossed-module $[F_2 \to
  F_1]$ by setting $F_2=F_1\times_{G_1} G_2$.  Then, we have a
  trivial fibration $[F_2 \to F_1] \to  [G_2 \to G_1]$. By
  assumption, this has a section, so in particular we get a section
  $s_1 \: G_1 \to F_1$ which embeds $G_1$ as a subgroup of $F_1$. It
  follows from Nielsen's  theorem that $G_1$ is free.
 \end{proof}

 \begin{rem}
  It is easy to see that a 2-group $\bbG$ is cofibrant in the
  Moerdijk-Svensson structure if an only if the inclusion $* \to
  \bbG$ is a cofibration. So the definition of cofibrant is the same
  in the pointed category.
  Also, in the pointed category, all 2-groupoids are fibrant.
 \end{rem}

 \begin{ex}{\label{E:replacement}}\indent\par
  \begin{itemize}

   \item[$\mathbf{1.}$] Let $\bbG=[G_2\to G_1]$ be an arbitrary crossed-module.
      Let $F_1 \to G_1$ be a surjective map from a free group $F_1$,
      and set $F_2:=F_1\times_{G_1}G_2$. Consider the crossed-module
      $\bbF=[F_2 \to F_1]$.  Then  $\bbF$ is cofibrant, and the natural
      map $\bbF \to \bbG$ is a trivial fibration (Definition \ref{D:fibrations}).
      In other words, $\bbF \to \bbG$ is a cofibrant replacement for $\bbG$.

   \item[$\mathbf{2.}$] Let $\Gamma$ be a group, and $F/R\cong \Gamma$
   be a presentation of $\Gamma$ as a quotient of a free group $F$.
   Then the map of crossed-modules $[R\to F] \to [1\to \Gamma]$
   is a cofibrant replacement for $\Gamma$.
  \end{itemize}
 \end{ex}

\begin{defn}{\label{D:mapping}}
   Let $\bbH$ and $\bbG$ be 2-groups (or crossed-modules). Choose a
   cofibrant replacement $\bbF \to \bbH$ for $\bbH$, as in Example
   \ref{E:replacement}.$\mathbf{1}$. The {\em derived mapping
   groupoid} $\Map(\bbH,\bbG)$ is defined to be $\homo(\bbF,\bbG)$,
   where $\homo$ is as in Definition \ref{D:homo}.
\end{defn}

Observe that  in the Moerdijk-Svensson structure all 2-groups are
automatically fibrant, so in the above definition we do not need a
fibrant replacement for $\bbG$.

The derived mapping groupoid $\Map(\bbH,\bbG)$ depends  on the
choice of the cofibrant replacement $\bbF \to \bbH$, but it is
unique up to an equivalence of groupoids (which is itself unique up
to transformation). Another way of thinking about the derived
mapping groupoid $\Map(\bbH,\bbG)$ is that it gives a model for the
groupoid of weak morphisms from $\bbH$ to $\bbG$ and pointed weak
transformations between them.\footnote{Since we are working in the
{\em pointed} category, the modification are trivial. This is due to
he fact that, whenever $X$ and $Y$ are pointed connected homotopy
2-types, the pointed mapping space $\Hombf_*(X,Y)$ is a 1-type.}

 \begin{thm}[\oldcite{Notes}, Proposition 8.1]{\label{T:weak}}
  Let $\bbG$ and $\bbH$ be 2-groups. Then, there is a
  natural (up to homotopy)   equivalence of groupoids
      $$\Map(\bbH,\bbG) \simeq \homw_*(\bbH,\bbG).$$
 \end{thm}

The fact that derived mapping groupoids are the correct models for
homotopy invariant mapping spaces is  justified  by the following.

 \begin{prop}{\label{P:mapping}}
  Let $\bbH$ and $\bbG$ be 2-groups (or crossed-modules). We have a
  natural bijection
      $$\pi_0\Map(\bbH,\bbG) \cong [\bbH,\bbG]_{\tgp}.$$
 \end{prop}

 \begin{proof}
  First let us remark that the proposition is not totally obvious
  because the category of 2-groupoids with the Moerdijk-Svensson
  model structure  is {\em not} a monoidal model category.

  By Proposition \ref{P:htpyequivalent}, we have
    $$\pi_0\Map(\bbH,\bbG)\cong\pi_0\Hombf_*(N\bbH,N\bbG)
      \cong [N\bbH,N\bbG]_{\mathbf{SSet}_*}.$$
  By Proposition \ref{P:equiv},
  $[N\bbH,N\bbG]_{\mathbf{SSet}_*}\cong [\bbH,\bbG]_{\tgp}$.
 \end{proof}

In $\S$\ref{S:Butterfly} we give a canonical explicit model for
$\Map(\bbH,\bbG)$, which is what we want.

\section{Some group theory}{\label{S:group}}

In this section we introduce some basic group theoretic lemmas which
will be used  in the proof of Theorem \ref{T:butterfly}. The results
in this section are mostly of technical nature.

\subsection{Generalized semi-direct products}{\label{SS:Generalized}}

We define a generalized notion of semi-direct product of groups, and
use that to introduce a pushout construction for crossed-modules.

Let $H$, $G$ and $K$ be  groups, each equipped with a right action
of $K$, the one on $K$ itself being conjugation. We denote all the
actions by $-^k$ (even the conjugation one). Assume we are given a
$K$-equivariant diagram
  $$\xymatrix@=12pt@M=10pt{
     H \ar[r]^p  \ar[d]_d  &   G   \\
          K            &       }$$
in which we require the compatibility condition $g^{d(h)}=g^{p(h)} \
\big(:=p(h)^{-1}gp(h)\big)$ is satisfied for every $h\in H$ and
$g\in G$.

 \begin{defn}{\label{D:crossed}}
    The {\em semi-direct product} $K\ltimes^H G$ of $K$ and $G$ along $H$
    is defined to be $K\ltimes G/N$, where
      $$N=\big\{\big(d(h)^{-1},p(h)\big), \ h \in H \big\}.$$
 \end{defn}

For this definition to make sense, we have to verify that $N$ is a
normal subgroup of $K\ltimes G$. This is left to the reader. Hint:
show that $G$ centralizes $N$, and  an element  $k \in K$ acts by
   $$\big(d(h)^{-1},p(h)\big)
             \mapsto \big(d(h^k)^{-1},p(h^k)\big).$$

There are natural group homomorphisms $p' \: K \to K\ltimes^H G$ and
$d' \: G \to K\ltimes^H G$, making the following diagram commute
     $$\xymatrix@=12pt@M=10pt{
                 H \ar[r]^p\ar[d]_d  &   G \ar[d]^{d'}  \\
                 K \ar[r]_(0.4){p'}        &   K\ltimes^H G    }$$
There is also an action of $K\ltimes^H G$ on $G$ which makes the
above diagram equivariant. An element $(k,g) \in K\ltimes^H G$ acts
on $x \in G$ by sending it to $g^{-1}x^kg$. Indeed, $[d'\: G \to
K\ltimes^H G]$ is a crossed-module.

 \begin{lem}{\label{L:crossed}}
    In the above square, the induced map $\Coker(d) \to \Coker(d')$
    is an isomorphism and the induced map $\Ker(d) \to \Ker(d')$
    is surjective. The kernel of the latter is equal to
    $\Ker(p) \cap \Ker(d)$.
 \end{lem}

 \begin{proof}
   Straightforward.
 \end{proof}

The relative semi-direct product construction satisfies the obvious
universal property. Namely, to give a homomorphism $K\ltimes^H G \to
T$ to an arbitrary group $T$ is equivalent to giving a pair of
homomorphisms $\delta \: G \to T$ and $\varpi \: K \to T$ such that
 \begin{itemize}
   \item \hspace{0.2in} $\xymatrix@=12pt@M=10pt{
            H \ar[r]^p \ar[d]_d \ar@{}[rd]|-{\circlearrowright} &
                                      G\ar[d]^{\delta}  \\
            K\ar[r]_{\varpi}        &   T    }$

   \item $\delta(g^k)=\delta(g)^{\varpi(k)}$, for every $g \in G$
     and $k \in K$.

 \end{itemize}

 \begin{rem}{\label{R:I}}
  Consider the the subgroup $I=p\big(\Ker(d)\big)=\Ker{d'} \subseteq
  G$.  It is   a $K$-invariant central subgroup of $G$. If in the
  relative semi-direct product construction we replace $G$ by
  $G/I$ the outcome will be the same.
 \end{rem}

In the following lemma we slightly modify the notation and denote
$K\ltimes^H G$ by $K\ltimes^{H,p} G$.

 \begin{lem}{\label{L:isom}}
  Notation being as above, let $\theta \: K \to G$ be a crossed homomorphism
  (i.e. $\theta(kk')=\theta(k)^{k'}\theta(k)$). Consider
  the group homomorphism $q \: H \to G$ defined by
  $q(h)=p(h)\theta(d(h))$, and use it to form $K\ltimes^{H,q} G$.
  Also consider the new action
  of $K$ on $G$ given by $g^{*k}:=\theta(k)^{-1}g^{k}\theta(k)$.
  (The $*$ is used just to differentiate the new action from the old one.)
  Then, there is a natural isomorphism
    $$\begin{array}{rccc}
           \theta_* \: &K\ltimes^{H,q} G & \risom & K\ltimes^{H,p} G \\
           &     (k,g) &\mapsto& \big(k,\theta(k)g\big).
      \end{array} $$
  The map $\theta_*$ makes the following triangle commute:
    $$\xymatrix@R=12pt@C=-2pt@M=7pt{
             & G  \ar[dl]_{d'_q} \ar[dr]^{d'_p}&    \\
          K\ltimes^{H,q} G    \ar[rr]_{\theta_*}  &  &
                                           K\ltimes^{H,p} G   }$$
  Furthermore, we have the following commutative
  triangle of isomorphisms
    $$\xymatrix@R=12pt@C=-6pt@M=7pt{
          \Coker(d'_q)  \ar[rr]^{\sim} \ar[dr]_{\sim}
                     &  &  \ar[dl]^{\sim} \Coker(d'_p) \\
             & \Coker{d} &    }$$
  where the top row is induced by  $\theta_*$.
 \end{lem}

 \begin{proof}
   We use the universal property of the relative semi-direct
   product. To give a map from  $K\ltimes^{H,p} G$ to a group $T$ is
   equivalent to giving a pair $(\delta,\varpi)$ of maps $\delta
   \: G \to T$ and $\varpi \: K \to T$ satisfying the two conditions
   described in the paragraph just before the lemma. To such a pair,
   we can associate a new pair $(\delta',\varpi')$, with
   $\delta':=\delta$ and
   $\varpi'(k):=\varpi(k)\delta(\theta(k))$. It is easy to see that
   the pair $(\delta',\varpi')$ satisfies the two conditions
   required by the universal property of $K\ltimes^{H,q} G$.
   Similarly, we can go backwards from a pair $(\delta',\varpi')$
   for $K\ltimes^{H,q} G$ to a pair $(\delta,\varpi)$  for
   $K\ltimes^{H,p} G$. It is easy to see that this correspondence is
   realized by $\theta_*$. This proves that $\theta_*$ is an
   isomorphism.

   Commutativity of the triangles is obvious.  (Also see Lemma
   \ref{L:crossed}.)
 \end{proof}

We have the following converse for Lemma \ref{L:isom}

 \begin{lem}{\label{L:triangle}}
   Consider  the semi-direct product diagrams
    $$\xymatrix@=12pt@M=10pt{
       H \ar[r]^p    \ar[d]_d  &   G   \\
          K            &       } \ \ \ \ \ \ \ \
    \xymatrix@=12pt@M=10pt{
     H \ar[r]^q    \ar[d]_d  &   G   \\
          K            &       }$$
  where the $K$-actions on $H$ are the same.
  Assume we are given an isomorphism of groups
  $\vartheta \: K\ltimes^{H,q} G  \risom K\ltimes^{H,p} G $
  such that the triangles of Lemma \ref{L:isom}   commute.
  Denote $p(\Ker d)\subseteq G$
  by $I$.
  \begin{itemize}
   \item[$\mathbf{i.}$] If  $\Ker d \subseteq \Ker p$,
     then there is a unique crossed homomorphism
     $\theta \: K \to G$ such that $\vartheta=\theta_*$ (see  Lemma \ref{L:isom},
     and Remark \ref{R:I}).

   \item[$\mathbf{ii.}$] If $K$ is a free group, then there exists a (not
     necessarily unique) crossed homomorphism
     $\theta \: K \to G$ such that $\vartheta=\theta_*$.

   \item[$\mathbf{iii.}$] If for crossed-homomorphisms
     $\theta$ and $\theta'$ we have $\theta_*=\theta'_*$, then
     the difference of $\theta$ and $\theta'$ factors through
     $I\subseteq G$. That is,
     for every $k \in K$, $\theta^{-1}(k)\theta'(k)$ lies in
     $I)$.
  \end{itemize}
 \end{lem}

 \begin{proof}[Proof of] ($\mathbf{i}$).
     Pick an element $(k,g) \in K\ltimes^{H,q} G$, and let
     $\theta(k,1)=(k',g')$. (Note that $(k,1)$ and $(k',g')$ are
     just representatives for actual elements in the corresponding
     relative semi-direct product groups). By the commutativity of
     the above triangle, the images of $k$ and $k'$ are the same in
     $\Coker(d)$; that is, there exists $h \in H$   such that
     $kd(h)=k'$. So, after adjusting $(k',g')$ by the
     $\big(d(h)^{-1},p(h)\big) \in N$ (see Definition
     \ref{D:crossed}), we may assume $k'=k$; that is
     $\vartheta(k,1)=(k,g')$. Define $\theta(k)$ to be $g'$. It is
     easily verified that $\theta$ is a crossed homomorphism, and
     that $\theta_*=\vartheta$.

  \vspace{0.1in}

  \noindent {\em Proof of} ($\mathbf{ii}$).
   Replace $G$ by $G/I$ and apply  ($\mathbf{i}$) to obtain
   $\theta \: K \to G/I$. Then use
   freeness of $K$ to lift $\theta$ to $G$.
  \vspace{0.1in}

  \noindent {\em Proof of} ($\mathbf{iii}$). Easy.
 \end{proof}

\subsection{A pushout construction for crossed-modules}

Continuing with the set-up of the previous section, we now bring
crossed-modules into the picture. Namely, we assume that $[d\: H \to
K]$ is a crossed-module. (Note that the condition $\mathbf{CM2}$ of
crossed-modules ($\S$\ref{SS:CrossedMod}) is already part of the
hypothesis.) To be compatible with our crossed-module notation, let
us denote $H$, $K$, $G$ and $d$ by $H_2$, $H_1$, $G_2$ and ${}_{-}$
respectively. Recall that $[G_2 \to H_1\ltimes^{H2} G_2]$ is again a
crossed-module.

 \begin{defn}{\label{D:pushout}} Let $\bbH=[H_2 \to H_1]$
  be a crossed-module. Let $G_2$ be a group with an action of $H_1$,
  and $p \: H_2 \to G_2$   an $H_1$-equivariant group homomorphism
  such that for every $\beta \in H_2$ and $\alpha \in G_2$ we have
  $\alpha^{\underline{\beta}}=\alpha^{p(\beta)}$.
  We call the crossed-module $[G_2 \to H_1\ltimes^{H2} G_2]$ the
  {\em pushout} of $\bbH$ along $p$ and denote it  by $p_*\bbH$.
 \end{defn}

 \begin{lem}{\label{L:pushout}}
  There is a natural induced map of crossed-modules $p_{\diamond} \:
  \bbH \to p_*\bbH$. Furthermore, $\pi_1(p_{\diamond})$ is an
  isomorphism and $\pi_2(p_{\diamond})$ is surjective. The kernel of
  $\pi_2(p_{\diamond})$ is equal to $$\{\beta \in H_2 \ | \
  \underline{\beta}=1, \ p(\beta)=1\}.$$ In particular, if $p \: H_2
  \to G_2$ is injective, then $p_{\diamond} \: [H_2 \to H_1] \to
  [G_2 \to H_1\ltimes^{H2} G_2]$ is
  an equivalence of crossed-modules.
 \end{lem}

 \begin{proof}
   Straightforward.
 \end{proof}

Assume now that we are given two crossed-modules $\bbG=[G_2 \to
G_1]$, $\bbH=[H_2 \to H_1]$ and  a morphism $P \: \bbH \to \bbG$
between them.  This gives us a diagram
   $$\xymatrix@=12pt@M=10pt{
     H_2\ar[r]^{p_2} \ar[d]_{-}  &  G_2   \\
           H_1          &       }$$
like the one in the beginning of this section. We also have an
action of $H_1$ on $G_2$ with respect to which $p_2$ is
$H_1$-equivariant. Namely, for $\alpha \in G_2$ and $h \in H_1$, we
define $\alpha^h$ to be $\alpha^{p_1(h)}$, the latter being the
action in $\bbG$. So, we can form the crossed-module $[G_2 \to
H_1\ltimes^{H2} G_2]$.

Define the map $\rho \:  H_1\ltimes^{H2} G_2 \to G_1$ by
$\rho(h,\alpha):=p_1(h)\underline{\alpha}$. It is easily seen to be
well-defined. We obtain the following commutative diagram of
crossed-modules:  \label{diagram1}

     $$\xymatrix@C=12pt@R=4pt@M=6pt{ H_2
        \ar[dd]\ar[rd]\ar@/^/[rr]^(0.35){p_2} & & G_2 \ar[dd]\\
                         &   G_2 \ar[dd] \ar[ru]_(0.48){=} & \\
          H_1  \ar[rd]\ar@/^/[rr]_(0.35){p_1} |!{[ru];[rd]}\hole
                                                  & &  G_1   \\
                    & H_1\ltimes^{H2} G_2 \ar[ru]_(0.55){\rho}& }$$
Note that the front-left square is almost an equivalence of
crossed-modules (Lemma \ref{L:pushout}); it is an actual equivalence
if and only if $\pi_2P \: \pi_2\bbH \to \pi_2\bbG$ is injective. If
this is the case, the above diagram means that, up to equivalence,
we have managed to replace our crossed-module map $P \: \bbH \to
\bbG$ with one, i.e. $P_*\bbH \to \bbG$, in which $p_2$ is the
identity map (the front-right square).

\vspace{0.1in} \noindent{\em Notation.} If $P\:\bbH \to \bbG$ is a
map of crossed-modules, we use the notation $P_*\bbH$ instead of
$p_{2,*}\bbH$.

\vspace{0.1in} The next thing we consider is, how the pushout
construction for crossed-modules behaves with respect to pointed
transformations between maps.

 \begin{lem}{\label{L:transformation}}
    Let $P, Q \: \bbH \to \bbG$ be maps of crossed-modules and
    $\theta \: Q \Rightarrow P$ a pointed transformation between
    them (Definition \ref{D:transformation}). Then, we have the
    following commutative  diagram of maps of crossed-modules:
      $$\xymatrix@R=-4pt@C=0pt@M=6pt{ H_2 \ar[rrrr]^{q_2}\ar[dddd]
        &&&& G_2\ar[ddr]_{=} \ar[dddd]|!{[ddl];[ddrrr]}\hole\ar[drrr]^{=}&&&\\
                                                        &&&&&&& G_2\ar[dddd]\\
                  &&& H_2\ar[rr]^{p_2}\ar[dddd]&& G_2\ar[dddd]\ar[urr]_{=}&&\\
                                                                     &&&&&&&\\
                        H_1 \ar[rrrr]_(0.35){q_1} |!{[urrr];[ddrrr]}\hole&&&&
                        H_1\ltimes^{H_2,Q} G_2\ar[ddr]_(0.45){\theta_*}
                    |-{\cong} \ar[drrr]^(0.42){\rho_Q}|!{[ur];[ddr]}\hole&&&\\
                                                                 &&&&&&& G_1\\
             &&& H_1\ar[rr]_(0.4){p_1}&&H_1\ltimes^{H_2,P} G_2
                                               \ar[rru]_(0.65){\rho_P} &&   }$$
    in which the front faces compose to $P$ and the back faces compose to $Q$.
    Here $\theta_*$ is obtained by the construction of Lemma \ref{L:isom}
    applied to  $\theta \: H_1 \to G_2$. Furthermore, the
    following triangle commutes:
       $$\xymatrix@R=12pt@C=-6pt@M=7pt{
          \pi_1(Q_*\bbH)  \ar[rr]_{\sim}^{\pi_1(\theta_*)} \ar[dr]_{\sim}
                     &  &  \ar[dl]^{\sim} \pi_1(P_*\bbH)  \\
                               & \pi_1(\bbH)  &    }$$
 \end{lem}

 \begin{proof}
  This is basically a restatement of Lemma \ref{L:isom}. Only proof
  of the equality $\rho_Q=\rho_P\circ\theta_*$ is missing. To prove
  this, pick $(h,\alpha) \in H_1\ltimes^{H_2,Q} G_2$. Since
  $\theta_*(h,\alpha)=\big(h,\theta(h)\alpha\big)$, we have
  \begin{eqnarray}
   \rho_P(\theta_*\big(h,\alpha)\big)=p_1(h)\underline{\theta(h)}\underline{\alpha}
    =q_1(h)\underline{\alpha}=\rho_Q(h,\alpha). \nonumber
  \end{eqnarray}
 \end{proof}

 \begin{lem}{\label{L:converse}}
  Consider the commutative diagrams of Lemma \ref{L:transformation},
  but with $\vartheta$ instead of $\theta_*$. Assume $\pi_2 P \:
  \pi_2\bbH \to \pi_2\bbG$ is the zero homomorphism. Then, there is a
  unique transformation $\theta \: Q \Rightarrow P$
  such that $\vartheta=\theta_*$.
 \end{lem}

 \begin{proof}
  This is more or less a restatement of Lemma
  \ref{L:triangle}.$\mathbf{i}$, with $H_1$, $H_2$ and $G_2$ playing
  the roles of $K$, $H$ and $G$, respectively. More explicitly,
  construct $\theta \: H_1 \to G_2$ as in Lemma \ref{L:triangle}. It
  automatically satisfies condition $\mathbf{T2}$ of Definition
  \ref{D:transformation}. The fact that it satisfies $\mathbf{T1}$
  follows from the  definition of $\theta_*$; see Lemma \ref{L:isom}.
 \end{proof}

\section{Butterflies as weak morphisms}{\label{S:Butterfly}}

In this section, we give a description for the groupoid of weak maps
between two crossed-modules (Theorem \ref{T:butterfly}). The key is
the following definition.

 \begin{defn}{\label{D:butterfly}}
  Let $\bbG=[G_2\to G_1]$ and $\bbH=[H_2\to H_1]$
  be crossed-modules. By a {\em butterfly} from $\bbH$ to $\bbG$ we
  mean a commutative diagram of groups
   $$\xymatrix@C=8pt@R=6pt@M=6pt{ H_2 \ar[rd]^{\kappa} \ar[dd]
                          & & G_2 \ar[ld]_{\iota} \ar[dd] \\
                            & E \ar[ld]^{\sigma} \ar[rd]_{\rho}  & \\
                                       H_1 & & G_1       }$$
  in which both diagonal sequences are complexes, and the NE-SW
  sequence, that is, $G_2 \to E \to H_1$, is short exact. We require
  $\rho$ and $\sigma$ satisfy the following compatibility with
  actions. For every $x \in E$, $\alpha \in G_2$, and $\beta \in
  H_2$,
    $$\iota(\alpha^{\rho(x)})=x^{-1}\iota(\alpha) x, \ \
      \kappa(\beta^{\sigma(x)})=x^{-1}\kappa(\beta) x.$$
  We denote the above butterfly by the tuple
  $(E,\rho,\sigma,\iota,\kappa)$. A {\em morphism} between two
  butterflies $(E,\rho,\sigma,\iota,\kappa)$ and
  $(E',\rho',\sigma',\iota',\kappa')$ is an isomorphism $f \: E \to
  E'$ commuting with all four maps. We define $\M(\bbH,\bbG)$ to be
  the groupoid of butterflies from $\bbH$ to $\bbG$.
 \end{defn}

 \begin{rem}
  Our notion of butterfly should not be confused with J. Pradines'.
  The latter correspond to  morphisms in the localized category of
  topological groupoids obtained by inverting Morita equivalences.
 \end{rem}

 \begin{lem}{\label{L:commute}}
  In a butterfly $(E,\rho,\sigma,\iota,\kappa)$, every element in the
  image of $\iota$ commutes with every element in $\Ker\rho$.
  Similarly, every element in the image of $\kappa$ commutes with
  every element in $\Ker\sigma$. In particular, the elements in the
  images of $\iota$ commute with the elements in the image of
  $\kappa$.
 \end{lem}

 \begin{proof}
  Easy.
 \end{proof}

The following theorem explains why we  butterflies are interesting
objects: they correspond to weak morphisms between crossed-modules
($\S$\ref{S:Weak}).

 \begin{thm}{\label{T:butterfly}}
  There is an equivalence of groupoids, natural up to
  a natural homotopy,
      $$\Omega \:\Map(\bbH,\bbG) \to  \M(\bbH,\bbG).$$
 \end{thm}

 \begin{proof}
  We begin with the following observation. Consider the product
  crossed-module $[H_2\times G_2 \to H_1\times G_1]$. Then, to give
  a butterfly from $\bbH$ to $ \bbG$ is equivalent  to giving a
  triangle
    $$\xymatrix@C=-4pt@R=10pt@M=10pt{
            & H_2\times G_2\ar[ld]\ar[rd]^{(\partial_{\bbH},\partial_{\bbG})}  &   \\
             E \ar[rr]_(0.45){(\sigma,\rho)}    &   &   H_1\times G_1   }$$
  which  has the property that $\rho$ intertwines the action on
  $H_2\times G_2$ of $E$ (by conjugation) with the action of
  $H_1\times G_1$ (from the crossed-module structure) and that the
  projection map $[H_2\times G_2 \to E] \to [H_2 \to H_1]$ is an
  equivalence of crossed-modules. (For this we use Lemma
  \ref{L:commute}.)

  Let us now construct the functor $\Omega \:\Map(\bbH,\bbG) \to
  \M(\bbH,\bbG)$. Choose a cofibrant replacement $\bbF=[F_2 \to
  F_1]$ for $\bbH$. Then, $\Map(\bbH,\bbG)\cong \homo(\bbF,\bbG)$.
  So, we have to construct $\homo(\bbF,\bbG) \to \M(\bbH,\bbG)$. We
  use the pushout construction of $\S$\ref{SS:Generalized}. Namely,
  given a morphisms $P \: \bbF \to \bbG$, we pushout $\bbF$ along
  the map $F_2 \to H_2\times G_2$ to obtain the crossed module $H_2
  \times G_2 \to E$. This is defined to be $\Omega(P)$.  This way we
  obtain a   triangle
       $$\xymatrix@C=-4pt@R=10pt@M=10pt{
            & H_2\times G_2\ar[ld]\ar[rd]^{(\partial_{\bbH},\partial_{\bbG})}    &   \\
             E \ar[rr]_(0.45){\rho}     &     &   H_1\times G_1
             }$$
  Lemma \ref{L:pushout} implies that the induced map $\bbF \to
  [H_2\times G_2 \to E]$ is a weak equivalence of crossed modules.
  Therefore, $[H_2\times G_2 \to E] \to [H_2 \to H_1]$ is also an
  equivalence of crossed-modules. Thus, we see that
  $\Omega(P):=[H_2\times G_2 \to E]$ is really a butterfly from
  $\bbH$ to $\bbG$.  This defines the effect of $\Omega$ on objects.
  The effect on morphisms is defined in the obvious way
  (see Lemma \ref{L:transformation}).

  We need to show that $\Omega$ is essentially surjective, full and
  faithful. Essential surjectivity follows from the fact that $\bbF$
  is cofibrant, and fullness follows from Lemma
  \ref{L:triangle}.$\mathbf{ii}$. Let us prove the faithfulness. Let
  $\theta,\theta' \: F_1 \to G_2$ be transformations between $P,Q \:
  \bbF \to \bbG$ such that $\Omega(\theta)=\Omega(\theta')$. Define
  $\hat{\theta} \: F_1 \to H_2 \times G_2$   by
  $\hat{\theta}=(1,\theta)$. Define $\hat{\theta'} \: F_1 \to H_2
  \times G_2$ in the similar way. By hypothesis, we have
  $\hat{\theta}_*=\hat{\theta'}_* : E_Q \to E_P$. By Lemma
  \ref{L:triangle}.$\mathbf{iii}$, for every $x \in F_1$, the
  element $\hat{\theta}^{-1}(x)\hat{\theta'}(x)$ lies in  the image
  of $\pi_2\bbF$ under the map $F_2 \to H_2\times G_2$. Since $\bbF
  \to \bbH$ induces an isomorphism on $\pi_2$, the only element in
  this image which is of the form $(1,\alpha)$ is $(1,1)$. On the
  other hand, every element $\hat{\theta}^{-1}(x)\hat{\theta'}(x)$
  is of the form $(1,\alpha)$. We conclude that
  $\hat{\theta}^{-1}(x)\hat{\theta'}(x)=(1,1)$. Therefore,
  $\theta(x)=\theta'(x)$. This completes the proof.
 \end{proof}

 \begin{rem}
  We can think of the category of crossed-modules and weak maps as
  the localization of $\tgp$ with respect to equivalences of
  crossed-modules.  Theorem \ref{T:butterfly} says
  that every weak map $f \: \bbH \to \bbG$ can be
  canonically written as a fraction

         $$\xymatrix@C=16pt@R=4pt@M=8pt{ & \mathfrak{K} \ar[ld]_{\sim}
              \ar[rd] & \\     \bbH  \ar@{..>} [rr]_f && \bbG}$$
  where
    $$\mathfrak{K}:=\sigma^*\bbH=[E\ltimes_{H_1}H_2 \llra{pr_1} E]
                     \cong[H_2\times G_2 \llra{(\kappa,\iota)} E].$$
 \end{rem}

 \begin{ex}
  Let $X$ be a topological space, and $A$ an abelian
  group. We can use the above theorem to give a description of the
  second cohomology $H^2(X,A)$. Recall that this cohomology group is
  in bijection with the set of homotopy classes of maps $X \to
  K(A,2)$, where $K(A,2)$ stands for the Eilenberg-MacLane space.
  Since $K(A,2)$ is simply connected, we can, equivalently, work with
  pointed homotopy classes. Assume the 2-type of $X$ is represented by
  the crossed-module $\bbH=[H_2 \to H_1]$. Then, there is a natural
  bijection $H^2(X,A)\cong [\bbH,A]_{\tgp}$, where we think of $A$ as
  the crossed module $[A \to 1]$. From Theorem \ref{T:butterfly} we
  conclude that $H^2(X,A)$ is in natural bijection with the set of
  isomorphism classes of pairs $(E,\kappa)$, where $E$ is a central
  extension of $H_1$ by $A$, and $\kappa \: H_2 \to E$ is an
  $E$-equivariant homomorphism. Here, $E$ acts on itself by
  conjugation and on $H_2$ via the projection $E \to H_1$.

  Of course, we can take any crossed-module $\bbG$ instead of $K(A,2)$
  and this way we find a description of the non-abelian cohomology
  $H^1(X,\bbG)$.
 \end{ex}

\subsection{Cohomological point of view}{\label{SS:CohButterfly}}

We quote the following cohomological description of the set of
pointed homotopy classes of weak maps from $\bbH$ to $\bbG$, that
is, $[\bbH,\bbG]_{\cm}$ from \cite{El} (which  apparently goes back
to Joyal and Street; see [ibid.], $\S$3.5).

 \begin{thm}[\oldcite{El}, Theorem 2.7]{\label{T:Elgueta}}
  Let $\bbG$ and $\bbH$ be crossed-modules, and assume they are
  represented by triples $(\pi_1\bbG,\pi_2\bbG,[\alpha])$ and
  $(\pi_1\bbH,\pi_2\bbH,[\beta])$, respectively, where $\alpha$ and
  $\beta$ are 3-cocycles representing the corresponding Postnikov
  invariants (see $\S$\ref{S:Homotopical}). Then, there is a
  bijection between the set of pointed homotopy classes of weak maps
  from $\bbH$ to $\bbG$, that is, $[\bbH,\bbG]_{\cm}$, and the set of
  triples $(\chi,\lambda,[c])$, where $\chi \: \pi_1\bbH \to
  \pi_1\bbG$ is a group homomorphism, $\lambda \: \pi_2\bbH \to
  \pi_2\bbG$ is a $\chi$-equivariant homomorphism such that
  $[\chi^*(\alpha)]=[\lambda_*(\beta)]$ in
  $H^3(\pi_1\bbH,\pi_2\bbG^{\chi})$, and $[c]$ is the class, modulo
  coboundary, of a 2-cochain on  $\pi_1\bbH$ with values in
  $\pi_2\bbG^{\chi}$, such that $\partial c= \chi^*(\alpha) -
  \lambda_*(\beta)$. Here, $\pi_2\bbG^{\chi}$ stands for $\pi_2\bbG$
  made into a $\pi_1\bbH$-module via $\chi$.
 \end{thm}

We will say more about this in the case where $\bbH=\Gamma$ is an
ordinary group in sections \ref{S:Cohomological}-\ref{S:Compatible}.

\section{Basic facts about butterflies}{\label{S:More}}

In this section we investigate butterflies in more detail.

\subsection{Butterflies vs.\! spans}{\label{SS:spans}}

Theorem \ref{T:butterfly} can be interpreted as saying that once we
invert strict equivalences of crossed-modules in $\cm$, the
morphisms of the resulting localized category can be presented, in a
canonical way, by butterflies. In fact, every butterfly
$(E,\rho,\sigma,\iota,\kappa)$ gives rise to a span of strict
crossed-module morphisms
 $$\xymatrix@C=8pt@R=6pt@M=6pt{& H_2\times G_2 \ar@{} |-{\cong}
         [dddl] \ar[ld]_{pr_1} \ar[rd]^{pr_2} \ar[dd]^{\mu}& \\
                          H_2  \ar[dd] & & G_2  \ar[dd] \\
                      & E \ar[ld]^{\sigma} \ar[rd]_{\rho}  & \\
                                       H_1 & & G_1       }$$
The map $\mu \: H_2\times G_2 \to E$ is defined by $(\alpha,\beta)
\mapsto \kappa(\alpha)\iota(\beta)$. Note that, by Lemma
\ref{L:commute}, the images of $\iota$ and $\kappa$ in $E$ commute.
The componentwise action of $E$ on $H_2$ and $G_2$ makes
$\mathbb{E}:=[\mu \: H_2\times G_2 \to E]$ into a crossed-module. It
is easy to verify that $\mathbb{E} \to \bbH$ is an equivalence of
crossed-modules.

\subsection{Butterflies vs.\! cocycles}{\label{SS:cocycle}}

  In (\cite{Notes}, Definition 8.4), we give a definition of a weak
  morphism between crossed-modules $[H_2 \to H_1]$ and $[G_2\to
  G_1]$ in terms of certain cocycles. By definition, such a cocycle
  consists of a triple $(p_1,p_2,\varepsilon)$, where $p_1 \: H_1
  \to G_1$, $p_2 \: H_2 \to G_2$, and $\varepsilon \: H_1\times H_1
  \to G_2$ are pointed set maps satisfying certain axioms
  [ibid.].\footnote{There is an error in the set of axioms in the
  published version of [ibid.]. This is corrected in an erratum
  which is available on my webpage.} Given such a
  triple, one can recover the corresponding butterfly as follows.

  Let $E$ be the group that has $H_1\times G_2$ as
  the underlying  set and whose product is defined by
     $$(h,g)\cdot(h',g')
     :=\big(hh',\varepsilon_{h,h'}^{-1}g^{p_1(h')}g'\big).$$
  Define the group homomorphisms $\rho \: E \to G_1$ and  $\kappa \:H_2 \to E$
  by
      $$\rho(h,g)=p_1(h)\underline{g}, \ \  \kappa(h)=(\underline{h},p_2(h)^{-1}).$$
   The homomorphisms $\iota \: G_2 \to E$ and $\sigma \: E \to H_1$
   are the inclusion and the projection maps on the corresponding
   components.

  It is easy to verify that this gives rise to a butterfly, and that this
  construction takes a (pointed) transformation of weak maps of crossed-modules
  (see {\em loc.\;cit.}) to a morphism of butterflies.

  Conversely, given a butterfly $(E,\rho,\sigma,\iota,\kappa)$, any
  choice of a set theoretic section $s \: H_1 \to E$ for $\sigma \:
  E \to H_1$ gives rise to a cocycle $(p_1,p_2,\varepsilon)$ defined by:
   $$\begin{array}{l}
         p_1 \: H_1 \to G_1, \ \  h \mapsto \rho s(h),\\
    p_2 \: H_2 \to G_2, \ \ \alpha \mapsto
    \kappa(\alpha)^{-1}s(\underline{\alpha}), \\
     \varepsilon_{h,h'}=s(h')^{-1}s(h)^{-1}s(hh').
   \end{array}$$

  If we choose a different section $s'$ for $\sigma$, we find  a pointed
  transformation $\theta $ between $(p_1,p_2,\varepsilon)$ and
  $(p'_1,p'_2,\varepsilon')$ given by
    $$\theta \: H_1 \to G_2, \ \ h\mapsto s(h)^{-1}s'(h).$$

  We can summarize this discussion by saying that, to give a triple
  $(p_1,p_2,\varepsilon)$ as in (\cite{Notes}, Definition 8.4)
  is the same thing as to give a butterfly $(E,\rho,\sigma,\iota,\kappa)$
  together with a set theoretic splitting $s \: H_1 \to E$ of $\sigma$.
  In particular, we have an equivalence of groupoids
    $$\M(\bbH,\bbG)\cong \Hom_{Weak}(\bbH,\bbG)$$
  where the right hand side is the groupoid whose morphisms of
  triples  $(p_1,p_2,\varepsilon)$ and whose morphisms of pointed
  transformations $\theta$  as in {\em loc.\;cit}.

\subsection{The induced map on homotopy groups}{\label{SS:homotopy}}

By Theorem \ref{T:butterfly}, a butterfly
$\mathcal{P}=(E,\rho,\sigma,\iota,\kappa)$ from $\bbH$ to $\bbG$
induces a well-defined map $N\bbH \to N\bbG$ in the {\em homotopy
category} of simplicial sets, which should be thought of as the
``nerve of the weak map $\mathcal{P}$''. Indeed, any choice of a
set-theoretic section $s$ for the map $\sigma \: E \to H_1$ gives
rise to a natural simplicial map $N_s\mathcal{P} \: N\bbH \to
N\bbG$. Furthermore, if $s'$ is another choice of a section,  there
is a natural simplicial homotopy between $N_s\mathcal{P}$ and
$N_{s'}\mathcal{P}$.

In particular, a butterfly $\mathcal{P}$ induces natural
homomorphisms on  $\pi_1$ and $\pi_2$. We can in fact describe these
maps quite explicitly.

To define $\pi_1\mathcal{P} \: \pi_1\bbH \to \pi_1\bbG$, let $x$ be
an element in $\pi_1\bbH$ and choose $y \in E$ such that
$\sigma(y)=x$. Define $\pi_1\mathcal{P}(x)$ to be the class of
$\rho(y)$ in $\pi_1\bbG$. This is easily seen to be well-defined. To
define $\pi_2\mathcal{P} \: \pi_2\bbH \to \pi_2\bbG$, let $\beta$ be
an element in $\pi_2\bbH$ and consider $y=\kappa(\beta)\in E$. Since
$\sigma(y)=\partial(\beta)=1$, there exists a unique $\alpha \in
G_2$ such that $\iota(\alpha)=y$. We define
$\pi_2\mathcal{P}(\beta)$ to be $\alpha$. Note that
$\partial(\alpha)=\rho(y)=1$, so $\alpha$ is indeed in $\pi_2\bbG$.

\subsection{Homotopy fiber of a butterfly}{\label{SS:more}}

We saw in $\S$\ref{SS:homotopy} that a butterfly
$(E,\rho,\sigma,\iota,\kappa)$  from $\bbH$  to $\bbG$ induces a map
$N\bbH \to N\bbG$ of simplicial sets which is well-defined up to
simplicial homotopy. It is natural to ask what is the homotopy fiber
of this map. In this subsection we see that this homotopy fiber can
be recovered from the NW-SE sequence of the butterfly. We describe
two ways of doing this.

Let  $\mathcal{P}\: (E,\rho,\sigma,\iota,\kappa)$ be a butterfly
from $\bbH$  to $\bbG$. We define its {\em homotopy fiber}
$\mathfrak{F}=\mathfrak{F}_{\mathcal{P}}$ to be the following
2-groupoid. The set of objects of $\mathfrak{F}$ is $G_1$. Given $g,
g' \in G_1$, the set of 1-morphisms in $\mathfrak{F}$ from $g$ to
$g'$ is the set of all $x \in E$ such that $g\rho(x)=g'$. For every
two such $x,y \in E$, the set of 2-morphisms from $x$ to $y$ is the
set of all $\gamma \in H_2$ such that $x\kappa(\gamma)=y$. We depict
this 2-cell by
      $$\xymatrix@C=10pt@R=-4pt@M=6pt{ & \ar@{=>}[dd]^{\gamma}&\\
        g \ar@/^1.1pc/ [rr]^{x} \ar@/_1.1pc/ [rr]_{y} &&  g'.  \\
            & &    }$$
It is clear how to define compositions rules in $\mathfrak{F}$,
 except perhaps for horizontal composition of 2-morphisms. Consider
two 2-morphisms
    $$\xymatrix@C=10pt@R=-4pt@M=6pt{ & \ar@{=>}[dd]^{\gamma}&
                                         & \ar@{=>}[dd]^{\delta}& \\
        g \ar@/^1.1pc/ [rr]^{x} \ar@/_1.1pc/ [rr]_{y} &&
         g' \ar@/^1.1pc/ [rr]^{z} \ar@/_1.1pc/ [rr]_{t} &&  g''.  \\
            & &   && }$$
We define their composition to be
     $$\xymatrix@C=28pt@R=-4pt@M=6pt{ &
                          \ar@{=>}[dd]^{\gamma^{\sigma(z)}\delta}&\\
        g \ar@/^1.1pc/ [rr]^{xz} \ar@/_1.1pc/ [rr]_{yt} &&  g''.  \\
            & &    }$$
There is a natural (strict) morphism of 2-groupoids $\Phi \:
\mathfrak{F} \to \bbH$. To describe this map, we need to think of
$\bbH$ as a 2-groupoid with one object,  as in $\S$\ref{SS:Equiv}.
We recall how this works. The unique object of $\bbH$ is denoted
$\bullet$. The set of 1-morphisms of $\bbH$ is $H_1$. Given two
1-morphisms $h,h' \in H_1$, the set of 2-morphisms from $h$ to $h'$
is the set of all $\gamma \in H_2$ such that $h\partial(\gamma)=h'$.

The natural map $\Phi \: \mathfrak{F} \to \bbH$ is described as
follows:
    $$\xymatrix@C=10pt@R=-4pt@M=6pt{ & \ar@{=>}[dd]^{\gamma}& & &
                                           & \ar@{=>}[dd]^{\gamma}& \\
        g \ar@/^1.1pc/ [rr]^{x} \ar@/_1.1pc/ [rr]_{y} && g' &
        \mapsto &   \bullet \ar@/^1.1pc/ [rr]^{\sigma(x)}
                        \ar@/_1.1pc/ [rr]_{\sigma(y)} &&  \bullet.  \\
                                                       & &  && && }$$

 \begin{thm}{\label{T:htpyfiber}}
  The sequence $N\mathfrak{F} \llra{N\Phi} N\bbH \llra{N\mathcal{P}} N\bbG$,
  which is well-defined in
  the homotopy category of simplicial sets, is a homotopy fiber sequence.
 \end{thm}

In order to prove Theorem \ref{T:htpyfiber}, we recall a few facts
about homotopy fibers in $\tgpd$. Given a strict morphism $P \: \bbH
\to \bbG$ of 2-groupoids, and a  base point $\bullet$ in $\bbG$,
there is a standard model for the homotopy fiber of $P$ which is
given by the following strict fiber product:
    $$ \operatorname{Fib}_P:=*\times_{\bullet, \bbG,s}
                                        \bbG^I \times_{t,\bbG,P}
    \bbH.$$ Here $\bbG^I:=\homc(0 \to 1, \bbG)$ is the 2-groupoid
of 1-morphisms of $\bbG$, and $s,t \: \bbG^I \to \bbG$ are the
source and target functors. (For a more precise definition of
$\bbG^I$ see Definition \ref{D:homo2}.)

In the case where $\bbG$ and $\bbH$ are 2-groups associated to
crossed-modules $[G_2\to G_1]$ and $ [H_2\to H_1]$, the 2-groupoid
$\operatorname{Fib}_P$ is described more explicitly as follows:
 \begin{itemize}
   \item  Objects are elements of $G_1$.

   \item 1-morphisms from $g$ to $g'$ are pairs $(x,\alpha)
     \in H_1\times G_2$, as in the 2-cell
          $$
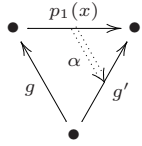
\xymatrix@=12pt{ \bullet   \ar[rr]^{p_1(x)} |
                                      (0.5){}="a"  &   & \bullet   \\
                                           &   &                   \\
        & \bullet \ar[luu]^g  \ar[ruu]_{g'} | (0.5){}="b"
                                          \ar @{:>}_{\alpha} "a";"b" &  }$$
     This means, $gp_1(x)\underline{\alpha}=g'$.

   \item 2-morphisms from $(x,\alpha)$ to $(y,\beta)$ are elements
    $\gamma \in H_2$ making the following diagram commutative:

        $$
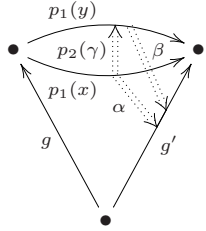
\xymatrix@=24pt{
                  \bullet    \ar@/_0.8pc/[rr]_(0.35){p_1(x)} | (0.55){}="d"
                \ar@/^0.8pc/[rr]^(0.35){p_1(y)} | (0.55){}="e"  | (0.63){}="b"
                                                          &   & \bullet   \\
                                                                   &  &   \\
            & \bullet \ar[luu]^g  \ar[uur]_{g'} | (0.65){}="a"  | (0.55){}="c"
               \ar @{:>}^(0.38){\beta} "b";"a" \ar @{:>}_(0.44){\alpha}"d";"c"
                                          \ar @{:>}^{p_2(\gamma)}  "d";"e"& }$$
      More precisely, we want $x\underline{\gamma}=y$ and
      $p_2(\gamma)\beta=\alpha$.
 \end{itemize}

There is a natural projection functor $pr
\:\operatorname{Fib}_P=*\times_{\bullet, \bbG,s}  \bbG^I
\times_{t,\bbG,f} \bbH \to \bbH$ which fits in the following fiber
homotopy sequence:
  $$\operatorname{Fib}_P \llra{pr} \bbH \llra{P} \bbG.$$

In the case where the map $p_2 \: H_2\to G_2$ is surjective, there
is a smaller model for the homotopy fiber of $P$ which we now
describe. Let $\operatorname{Fib}'_P$ be the full sub-2-category of
$\operatorname{Fib}_P$ in which for 1-morphisms we only take the
ones for which $\alpha$ is the identity. It is easily checked that,
in this case, the inclusion $\operatorname{Fib}'_P \subset
\operatorname{Fib}_P$ is an equivalence of 2-groupoids.

Let us record this as a lemma.

 \begin{lem}{\label{L:fiber}}
  Let $\bbG=[G_2\to G_1]$ and $\bbH =[H_2\to H_1]$ be crossed-modules,
  viewed as 2-groups.
  Let $P \: \bbH \to \bbG$ be a strict morphism of 2-groups,
  and let $\operatorname{Fib}'_P$ be the 2-groupoid defined in the previous
  paragraph. If  $p_2 \: H_2\to G_2$ is surjective, then the sequence
    $$\operatorname{Fib}'_P \llra{pr} \bbH \llra{P} \bbG$$
  is a homotopy fiber sequence
  in the category of 2-groupoids.
 \end{lem}

We can now prove Theorem \ref{T:htpyfiber}.

 \begin{proof}[Proof of Theorem \ref{T:htpyfiber}]
   We apply Lemma \ref{L:fiber} to the morphism  of crossed-modules
   $P \: [H_2\times G_2 \to E] \to [G_2 \to G_1]$. Observe that
   the  crossed-module on the left is a model for $\bbH$
   (via the natural weak equivalence
   $(pr_1, \sigma) \: [H_2\times G_2 \to E] \to [H_2\to H_1]$),
   so $P$ is a strict model for $\mathcal{P}$.

   It is easily verified that the 2-groupoid $\mathfrak{F}$ of the theorem
   is canonically {\em isomorphic} to the 2-groupoid
   $\operatorname{Fib}'_P$ of Lemma \ref{L:fiber}. So, the homotopy
   fiber sequence of Lemma \ref{L:fiber} is naturally isomorphic, in
   the homotopy category of 2-groupoids, to the sequence
   $$\mathfrak{F}
   \llra{\Phi} \bbH \llra{\mathcal{P}} \bbG.$$ Taking the nerves
   gives the homotopy fiber sequence  of Theorem \ref{T:htpyfiber}.
 \end{proof}

We now derive some corollaries of Theorem \ref{T:htpyfiber}.

 \begin{prop}{\label{P:butterflyfiber}}
  Let $\mathcal{P}\: (E,\rho,\sigma,\iota,\kappa)$ be a  butterfly
  from $\bbH$ to $\bbG$. Consider the nerve $N\mathcal{P} \: N\bbH
  \to N\bbG$, which is well-defined in  the homotopy category of
  simplicial sets,  and let $F$ be its homotopy fiber (note that $F$
  is naturally pointed). Let
    $$C: \ \  H_2 \llra{\kappa} E \llra{\rho} G_1$$
  be the NW-SE sequence of $\mathcal{P}$. Then, there are natural
  isomorphisms $H_i(C)\cong \pi_iF$, $i=0,1,2$.
  (Of course, for $i=0$ this means an isomorphism of pointed sets.)
 \end{prop}

 \begin{proof}
  By Theorem \ref{T:htpyfiber}, $F$ is naturally homotopy equivalent
  to $N\mathfrak{F}$.
  It is easily verified that $H_i(C)\cong \pi_i\mathfrak{F}$, $i=0,1,2$.
  The result follows.
 \end{proof}

 \begin{rem}
   Note that the truncated sequence $H_2 \to E$ is indeed naturally
   a crossed-module (take the  action of $E$ on $H_2$ induced via
   $\sigma$). The 2-groupoid $\mathfrak{F}$ of Theorem \ref{T:htpyfiber} can
   be recovered from the sequence $C$ together with this extra structure.
 \end{rem}

 \begin{cor}{\label{C:sequence}}
  Notation being as in Proposition \ref{P:butterflyfiber},
  there is a long exact sequence
   $$\xymatrix@=7pt{ 1 \ar[r] & H_2(C) \ar[r]
                 &\pi_2\mathfrak{H} \ar[r] &\pi_2\mathfrak{G} \ar[r]
                 &H_1(C) \ar[r] &\pi_1\mathfrak{H} \ar[r]
                 &\pi_1\mathfrak{G} \ar[r] & H_0(C) \ar[r] &1. }$$
 \end{cor}

 \begin{proof}
  Immediate.
 \end{proof}

 \begin{prop}{\label{P:butterflyequiv}}
  A butterfly $(E,\rho,\sigma,\iota,\kappa)$   from $\bbH$
  to $\bbG$ is a weak equivalence (i.e. induces a homotopy
  equivalence on the nerves), if and only if the NW-SE sequence
     $$   H_2 \llra{\kappa} E \llra{\rho} G_1$$
  is short exact. A weak inverse for this butterfly in obtained
  by simply flipping the butterfly, as in the following diagram
   $$\xymatrix@C=8pt@R=6pt@M=6pt{ G_2 \ar[rd]^{\iota} \ar[dd]
                            & & H_2 \ar[ld]_{\kappa} \ar[dd] \\
                  & E \ar[ld]^{\rho} \ar[rd]_{\sigma}  & \\
                  G_1 & & H_1       }$$
 \end{prop}

 \begin{proof}
  Immediate.
 \end{proof}

There is another way of describing the homotopy fiber of a
butterfly. Define $\bbF$ to be the following crossed-module in
groupoids:
   $$\bbF:=[\underset{g \in G_1}{\coprod}H_2
                        \llra{f} G_1\times E \sst{s,t} G_1].$$
The groupoid $[G_1\times E \sst{} G_1]$ is the translation groupoid
of the right multiplication action of $E$ on $G_1$ via $\rho$. Its
source and target maps are $s\: (g,x) \mapsto g$ and $t\: (g,x)
\mapsto g\rho(x)$. The restriction of  $f$ to the copy of $H_2$
corresponding to the index $g \in G_1$ sends the element $\beta \in
H_2$ to $\big(g,\kappa(\beta)\big)$.

It can be shown that the nerve of $\bbF$ is naturally equivalent to
$\mathfrak{F}$.

\subsection{Weak morphisms vs.\! strict morphisms}{\label{SS:st}}

We saw in the Theorem \ref{T:butterfly} that $\M(\bbH,\bbG)$ is a
model for the groupoid of weak morphism from $\bbH$ to $\bbG$.
Strict morphisms $\bbH \to \bbG$ form a  subgroupoid of
$\M(\bbH,\bbG)$. The objects of this  subgroupoid are precisely the
butterflies for which the NE-SW short exact sequence is split.

More precisely, given a strict morphisms $P \: \bbH \to \bbG$, the
corresponding butterfly looks as follows:
$$\xymatrix@C=4pt@R=2pt@M=10pt{ H_2 \ar[rd]^{\kappa} \ar[dd]
                            & & G_2 \ar[ld]_{\iota} \ar[dd]  \\
                  & H_1\ltimes G_2 \ar[ld]^{\sigma} \ar[rd]_{\rho}  & \\
                  H_1 & & G_1       }$$
where $\iota=(1,\id)$, $\sigma=pr_1$,
$\kappa(\beta)=\big(\underline{\beta},p_2(\beta^{-1})\big)$, and
$\rho(x,\alpha)=p_1(x)\underline{\alpha}$. Here, the action of $H_1$
on $G_2$ is obtained via $p_1$ from that of $G_1$ on $G_2$.

Observe that, any butterfly coming from a strict morphism has a
canonical splitting. A morphism of butterflies does not necessarily
respect this splitting. In fact, the difference between the
resulting splittings  determines a (pointed) transformation
(Definition \ref{D:transformation}) between the corresponding strict
functors, and vice versa. Let us state this in the following.

\begin{prop}{\label{P:splitting}}
 Let $\mathcal{SB}(\bbH,\bbG)$ be the groupoid whose objects are
 butterflies with a splitting and whose morphisms are the morphisms
 of the underlying butterflies. Then, we have a natural equivalence of
 groupoids
   $$\homo(\bbH,\bbG)\cong\mathcal{SB}(\bbH,\bbG).$$
\end{prop}


\subsection{Butterflies in the differentiable or algebraic
contexts}{\label{SS:Lie}}

Our discussion of butterflies can be generalized to a global setting
in which everything happens over a Grothendieck site $\mathsf{S}$. A
group is now replaced by a sheaf of groups, and a crossed-module
$[G_2\to G_1]$ will have $G_1$ and $G_2$ sheaves of groups. Also, in
the definition of a butterfly we require that the NW-SE sequence is
a short exact sequence of sheaves of groups.

It can be shown that the discussion of the previous, and the
subsequent, sections on butterflies goes through more or less
verbatim in the relative case. I particular, it can be shown that
butterflies model morphisms in the homotopy category of
crossed-modules over $\mathsf{S}$.

As a consequence of this global approach, we obtain  theories of
{\em butterflies for crossed-modules in Lie groups, topological
groups, group schemes} and so on. For example, a Lie butterfly
between two Lie crossed-modules is one in which $E$ is a Lie group,
all the maps $\iota$, $\kappa$, $\rho$, and $\sigma$ are
differentiable homomorphisms, and the sequence
    $$1 \to G_2 \llra{\iota} E \llra{\sigma} H_1 \to 1$$
is a short exact sequence of Lie groups.

The general theory of butterflies over a Grothendieck site will
appear in a joint paper with E.~Aldrovandi.

\begin{rem}
  The savvy reader may complain that in the example of Lie
  butterflies mentioned above, one can only expect $E$ to be
  a sheaf of groups over the site $\mathsf{S}$ of differentiable
  manifolds. However, it is not hard to see that when $E$ is a
  sheaf of groups which sits in a short exact sequence whose both
  ends are Lie groups, then $E$ itself is necessarily representable
  by a Lie groups.
\end{rem}

\section{Bicategory of crossed-modules and weak
maps}{\label{S:Bicat}}

In this section we construct a bicategory $\mathcal{CM}$ whose
objects are crossed-modules and whose 1-morphisms are butterflies.
The bicategory $\mathcal{CM}$ is a model for the homotopy category
of pointed connected 2-types.

\subsection{The bicategory $\mathcal{CM}$}{\label{SS:CM}}

Theorem \ref{T:butterfly} enables us to give a concrete model for
the 2-category of crossed-modules, weak morphisms and weak
transformation. More precisely, define $\mathcal{CM}$ to be the
bicategory whose objects are crossed-modules and whose morphism
groupoids are $\M(\bbH,\bbG)$. The composition functors
$\M(\bbK,\bbH)\times \M(\bbH,\bbG) \to \M(\bbK,\bbG)$ are defined as
follows. Given butterflies
  $$\xymatrix@C=8pt@R=6pt@M=6pt{ K_2 \ar[rd]  \ar[dd]_{\xi}
                            & & H_2 \ar[ld]_{\iota'}  \ar[dd]^{\partial} \\
                  & F \ar[ld]  \ar[rd]_{\rho'}   & \\
                  K_1 & & H_1       } \ \ \ \ \
 \xymatrix@C=8pt@R=6pt@M=6pt{ H_2 \ar[rd]^{\kappa}  \ar[dd]_{\partial}
                            & & G_2 \ar[ld]  \ar[dd]^{\partial} \\
                  & E \ar[ld]^{\sigma}  \ar[rd]  & \\
                  H_1 & & G_1       }$$
we define their composite to be the butterfly
$$\xymatrix@C=10pt@R=4pt@M=4pt{ K_2 \ar[rd]  \ar[dd]_{\xi}
                            & & G_2 \ar[ld]  \ar[dd]^{\partial} \\
                  & F\oux{H_1}{H_2}E \ar[ld]  \ar[rd]   & \\
                  K_1 & & G_1       }$$
where $F\oux{H_1}{H_2}E $ is a Baer type product. More precisely, it
is the quotient of the group
  $$L:= \{(y,x) \in F\times E \ | \ \rho'(y)=\sigma(x) \in H_1\}$$
modulo the subgroup
   $$I=\{\big(\iota'(\beta),\kappa(\beta)\big)
  \in F\times E\ | \ \beta \in H_2\}.$$

We have the following theorem (see $\S$\ref{S:Transfo} for
notation).

 \begin{thm}{\label{T:bicat}}
  With morphism groupoids being $\M(\bbH,\bbG)$,
  crossed-modules form a bicategory $\mathcal{CM}$.
  There is a natural weak functor $\underline{\cm} \to \mathcal{CM}$
  which is fully faithful on morphism
  groupoids.
 \end{thm}

 \begin{proof}
   It is straightforward to verify that $\mathcal{CM}$ is a
   bicategory. The functor $\underline{\cm} \to \mathcal{CM}$
   is the one constructed in  $\S$\ref{SS:st}. Fully faithfulness
   on morphism groupoids is a restatement of Proposition
   \ref{P:splitting}.
 \end{proof}

 \begin{rem}
  The bicategory $\mathcal{CM}$ is a model for (i.e., is naturally
  biequivalent to) the 2-category of 2-groups and {\em weak}
  morphisms between them; see \cite{Notes}, especially Propositions
  8.1 and  7.8. More precisely, the latter is the 2-category whose
  objects are 2-groups, whose 1-morphisms are weak morphisms between
  2-groups, and whose 2-morphisms are pointed
  transformations between them.

  The advantage of working with $\mathcal{CM}$ is that 1-morphisms
  and 2-morphisms in $\mathcal{CM}$ are quite explicit, and in order
  to describe weak morphism one does not need to deal
  with complicated cocycles and  coherence conditions (like the ones coming from
  the hexagon axiom) which in practice make
  explicit computations intractable. This is a great advantage in geometric
  contexts, say, when dealing with weak morphisms of Lie or algebraic 2-groups
  (see \ref{SS:Lie}).
\end{rem}

\subsection{A special case of composition of butterflies}
{\label{SS:specialcompose}}

There is a simpler, and quite useful, description for the
composition of two butterflies in the case where one of them is
strict ($\S$\ref{SS:st}).

When the first morphisms is strict, say
  $$\xymatrix@C=8pt@R=6pt@M=6pt{ K_{2}  \ar[rr]^{q_{2}} \ar[dd]
                              & & H_{2}   \ar[dd]  \\
                                           &     & \\
                        K_1 \ar[rr]_{q_1} & & H_1       }$$
then the composition is
      $$   \xymatrix@C=8pt@R=6pt@M=6pt{ K_{2} \ar[rd]   \ar[dd]
                            & & G_{2} \ar[ld]  \ar[dd]  \\
      & q^{*}_1(F) \ar[ld]^{q^{*}_1(\sigma')} \ar[rd] & \\
                         K_1 & & G_1       }$$
Here, $q^{*}_1(F)$ stands for the pull back of the extension $F$
along $q_1 \: K_1 \to H_1$. More precisely,
$q^{*}_1(F)=K_1\times_{H_1}F$ is the fiber product.

When the second morphisms is strict, say
          $$\xymatrix@C=8pt@R=6pt@M=6pt{ H_{2}  \ar[rr]^{p_{2}}
          \ar[dd]
                            & & G_{2}   \ar[dd]  \\
                  &          & \\
                  H_1 \ar[rr]_{p_1} & & G_1       }$$
 then the composition is
         $$   \xymatrix@C=8pt@R=6pt@M=6pt{ K_{2} \ar[rd]   \ar[dd]
                    & & G_{2} \ar[ld]_{p_{2,*}(\iota)}  \ar[dd]  \\
                                & p_{2,*}(E) \ar[ld]  \ar[rd]   & \\
                                              K_1 & & G_1      }$$
Here, $p_{2,*}(E)$ stands for the push forward of the extension $E$
along $p_{2} \: H_{2} \to G_{2}$. More precisely,
$p_{2,*}(E)=E\times_{H_{2}} G_{2}$ is the pushout.

\section{Kernels and cokernels of butterflies}{\label{S:Ker}}
%

In this section we give an overview of how certain basic notions of
group theory, such as kernels, cokernels, complexes, extensions, and
so on, can be generalized to the setting of 2-groups and weak
morphisms through the language of butterflies. These notions have,
of course, been studied by various authors already, but our emphasis
is on showing how our approach via butterflies makes things
remarkably simple and reveals structures that are not easy to see
using the cocycle approach.

This section is meant as a general overview and we do not give many
proofs. The proofs are not at all hard, but we believe this topic
deserves a systematic treatment which is out of the scope of this
paper. This is being worked out in \cite{Extensions}.

\subsection{Kernels and cokernels of butterflies}
{\label{SS:kernelbutterfly}}

For a weak morphism $\bbH \to \bbG$ of 2-groups, one can define a
kernel and a cokernel using a certain universal property. The kernel
is always expected to be a 2-group again. The cokernel, however, is
only a pointed groupoid in general. In this section, we give an
explicit description of the kernel and the cokernel of a weak
morphism of 2-groups using the formalism of butterflies.

Consider the butterfly $\mathcal{P}\: \bbH \to \bbG$ given by
  $$\xymatrix@C=8pt@R=6pt@M=6pt{ H_2 \ar[rd]^{\kappa} \ar[dd]_{\partial}
                            & & G_2 \ar[ld]_{\iota} \ar[dd]^{\partial} \\
                  & E \ar[ld]^{\sigma} \ar[rd]_{\rho}  & \\
                  H_1 & & G_1       }$$
We define the {\em kernel} of $\mathcal{P}$ by
  $$\Ker\mathcal{P}:=[H_2\llra{\kappa}\Ker\rho],$$
and the {\em cokernel} of $\mathcal{P}$ by
    $$\Coker\mathcal{P}:=[\Coker\kappa \llra{\rho} G_1].$$
Note that $\Coker\mathcal{P}$ is  just a group homomorphism and  may
not be a crossed-module  in general (also see $\S$\ref{SS:braided}).

\begin{rem}{\label{R:fiber}}
  The 2-group associated to the crossed-module $\Ker\mathcal{P}$ is
  naturally equivalent to the 2-group of automorphisms of the base
  point of the homotopy fiber $\mathfrak{F}_{\mathcal{P}}$ of
  $\mathcal{P}$ defined in $\S$\ref{SS:more}. The pointed set
  $\pi_1(\Coker\mathcal{P})$ is in natural bijection with the set of
  connected components of $\mathfrak{F}_{\mathcal{P}}$. We have
  a natural isomorphism
  $$\pi_1\Ker\mathcal{P}\cong\pi_2\Coker\mathcal{P}.$$
  (Note: using the notation $\pi_1$ and $\pi_2$ for the cokernel and the
  kernel of a group homomorphism $[K \to L]$ which is not a
  crossed-module (e.g., for $[K \to L]=\Coker\mathcal{P}$) is not quite
  appropriate. We have used it here for the sake of notational consistency.
  This will not appear elsewhere in the paper except for this subsection.)
 \end{rem}

There is a (strict) morphism $I_{\mathcal{P}}:=(\id_{H_2},\sigma) \:
\Ker\mathcal{P} \to  \bbH$.  The morphism $I_{\mathcal{P}}$ comes
with a natural 2-morphism $\theta_{\mathcal{P}} \: \mathcal{P}\circ
I_{\mathcal{P}} \Rightarrow 1$ in
$\Hom_{\mathcal{CM}}(\Ker\mathcal{P},\bbG)$, where $1$ stands for
the trivial morphism. The following proposition follows from
Proposition \ref{P:trivial} below. We leave it to the reader to
supply the proof and also to formulate the corresponding statement
for the cokernel of a butterfly.

 \begin{prop}{\label{P:kernel}}
  The pair $(I_{\mathcal{P}},\theta_{\mathcal{P}})$ satisfies the
  universal property of a kernel for
  $\mathcal{P} \: \bbH \to \bbG$ in the bicategory
  $\mathcal{CM}$. That is, for every crossed-module $\mathfrak{L}$,
  $(I_{\mathcal{P}},\theta_{\mathcal{P}})$
  induces an equivalence  of groupoids between
  $\Hom_{\mathcal{CM}}(\mathfrak{L},\Ker\mathcal{P})$
  and the homotopy fiber of $\mathcal{P}_* \: \Hom_{\mathcal{CM}}(\mathfrak{L},\bbH)
  \to \Hom_{\mathcal{CM}}(\mathfrak{L},\bbG)$ over the trivial map
      $1 \in  \Hom_{\mathcal{CM}}(\mathfrak{L},\bbG)$.
 \end{prop}

 \begin{prop}{\label{P:sequence}}
  Let $\mathcal{P} \: \mathfrak{H} \to \mathfrak{G}$ be a butterfly.
  The map $I_{\mathcal{P}}\: \Ker\mathcal{P} \to \bbH$
  gives rise to a long exact sequence:
   $$\xymatrix@=7pt{ 1 \ar[r] &\pi_2\Ker\mathcal{P} \ar[r]
                 &\pi_2\mathfrak{H} \ar[r] &\pi_2\mathfrak{G} \ar[r]
                 &\pi_1\Ker\mathcal{P} \ar[r] &\pi_1\mathfrak{H} \ar[r]
                 &\pi_1\mathfrak{G} \ar[r] &\pi_1\Coker\mathcal{P} \ar[r] &1. }$$
 \end{prop}

 \begin{proof}
  Exercise.
 \end{proof}

 \begin{cor}{\label{C:equivalence}}
     A butterfly $\mathcal{P} \: \mathfrak{H} \to \mathfrak{G}$ is
     an equivalence if and only if $\Ker\mathcal{P}$ and
     $\Coker\mathcal{P}$ are    trivial (i.e., equivalent to a point).
 \end{cor}


\begin{rem}
  One can also define the image and the coimage of a butterfly.
  Kernel, image,  cokernel, and the coimage of a butterfly correspond
  to the top-left, top-right, bottom-right, and the bottom-left
  arrows in the butterfly, respectively.

  In fact, careful study of butterflies
  shows that each of the notions kernel,  image, cokernel, and the
  coimage can be defined in {\em two} ways, one of which we have given above.
  So, we have eight
  different crossed-modules associated to a butterfly, and these
  eight crossed-modules interact in an interesting way.
  This seems to suggest that one ought to consider all eight at once. This is
  discussed in detail in \cite{Extensions}.
\end{rem}

\subsection{Exact sequences of crossed-modules in $\mathcal{CM}$}{\label{SS:exact}}

The following proposition is the key in studying complexes of
2-groups.

 \begin{prop}{\label{P:trivial}}
  Let $\mathcal{Q}$ and $\mathcal{P}$ be butterflies as in the
  following diagram:
  $$\xymatrix@C=8pt@R=6pt@M=6pt{ K_2 \ar[rd]^{\kappa'}  \ar[dd]
                            & & H_2 \ar[ld]_{\iota'}
                            \ar[rd]^{\kappa}  \ar[dd]
                                 && G_2 \ar[ld]  \ar[dd]  \\
                  & F \ar[ld] \ar[rd]_{\rho'}
                  \ar@{..>}[rr]^(0.38){\delta} |!{[ur];[dr]}\hole
                   & & E  \ar[ld]^{\sigma}   \ar[rd]_{\rho}  & \\
                  K_1 & & H_1    & & G_1   }$$
  Then, the composition $\mathcal{P}\circ\mathcal{Q}$ is the trivial
  butterfly if and only if there exists a group homomorphism $\delta
  \: F \to E$ making the above diagram commutative such that
  the sequence
      $$1 \llra{} K_2 \llra{\kappa'} F \llra{\delta} E
                                   \llra{\rho} G_1 \llra{} 1$$
  is a complex (i.e., the composition of consecutive maps is the
  trivial map). 
 \end{prop}

Let $\mathcal{Q} \: \bbK \to \bbH$ and  $\mathcal{P} \: \bbH \to
\bbG$ be butterflies. Proposition \ref{P:trivial} gives a necessary
and sufficient condition for the composition
$\mathcal{P}\circ\mathcal{Q}$ to be zero. Observe that when this
condition is satisfied, we have a natural  morphism $\bbK \to
\Ker\mathcal{P}$ given by the butterfly
   $$\xymatrix@C=8pt@R=6pt@M=6pt{ K_2 \ar[rd]^{\kappa'} \ar[dd]
                          & & H_2 \ar[ld]_{\iota'} \ar[dd]^{\kappa} \\
                            & F \ar[ld]^{\sigma'} \ar[rd]_{\delta}  & \\
                                       K_1 & & \Ker\rho      }$$

The next question to ask is when the sequence
   $$\bbK \llra{\mathcal{Q}} \bbH \llra{\mathcal{P}} \bbG$$
is exact at $\bbH$. The correct definition of exactness seems to be
the following.

\begin{defn}{\label{D:exact}}
  We say that the sequence
     $$\bbK \llra{\mathcal{Q}} \bbH \llra{\mathcal{P}} \bbG$$
  is {\bf exact} at $\bbH$ if $\mathcal{P}\circ\mathcal{Q}$ is zero and $\Coker(\bbK \to
\Ker\mathcal{P})$ is trivial.
\end{defn}

\begin{rem}
  If we pretend for a moment that
  $\Coker\mathcal{Q}=[\Coker{\kappa'} \llra{\rho'} H_1]$ is a
  crossed-module,  then we obtain a ``butterfly'' $\Coker\mathcal{Q} \to \bbG$
  given by
   $$\xymatrix@C=8pt@R=6pt@M=6pt{
                   \Coker{\kappa'} \ar[rd]^{\delta} \ar[dd]_{\rho'}
                          & & G_2 \ar[ld]_{\iota} \ar[dd] \\
                           & E \ar[ld]^{\sigma} \ar[rd]_{\rho}  & \\
                                       H_1 & & G_1       }$$
   Definition \ref{D:exact} is now equivalent to
   $\Ker(\Coker\mathcal{Q} \to \bbG)$ being trivial. This argument
   is actually a valid argument if we assume $\bbH$ is a braided
   crossed-module because in this case $\Coker\mathcal{Q}$
   is, indeed, a crossed-module and the above diagram is an honest butterfly;
   see $\S$\ref{S:Derived}.
\end{rem}

The following proposition is immediate from the definition.

\begin{prop}{\label{P:exact}}
  The sequence $\bbK \llra{\mathcal{Q}} \bbH \llra{\mathcal{P}}
  \bbG$ is exact at $\bbH$ (respectively, short exact) if and only
  if the sequence $$1\to  K_2 \llra{\kappa'} F \llra{\delta} E
  \llra{\rho} G_1 \to 1$$ of Proposition \ref{P:trivial} is exact at
  $F$ and $E$  (respectively, exact everywhere).
\end{prop}

Using the above proposition the reader can work out what a complex
(or an exact sequence) of crossed-modules in $\mathcal{CM}$ looks
like.

Proposition \ref{P:exact} also provides a new perspective on the
problem of studying extensions of a 2-group $\bbG$ by a 2-group
$\bbK$ (see \cite{Breen2} and \cite{Rou}). This will be studied in
more detail in \cite{Extensions}.

\section{Braided and abelian butterflies}{\label{S:Derived}}

In this section we will discuss butterflies in an abelian category
$\A$. This is intended to be an overview of the main features of
braided butterflies and we will not give proofs.

We obtain an explicit description of the derived category of
complexes of length two in $\A$; compare \cite{Del}. We begin by
discussing braided butterflies.

\subsection{Braided butterflies}{\label{SS:braided}}

Assume $\bbH$ and $\bbG$ are braided crossed-modules \cite{Co}, and
let $\{-,-\}_{\bbH}$ and $\{-,-\}_{\bbG}$ denote the corresponding
braidings. (Our braiding convention is that
$\partial\{x,y\}=xyx^{-1}y^{-1}$.)

 \begin{defn}{\label{D:braided}}
  A butterfly $\mathcal{P} =(E,\rho,\sigma,\iota,\kappa)\:\bbH
  \to\bbG$ is called {\bf braided} if the following conditions are
  satisfied:
    $$\forall x,y \in E, \ \
    \kappa\{\sigma(x),\sigma(y)\}_{\bbH}=
         xyx^{-1}\iota\{\rho(x),\rho(y^{-1})\}_{\bbG}y^{-1}.$$
 \end{defn}

In the case where $\mathcal{P}$ comes from a strict morphism of
2-groups ($\S$\ref{SS:st}) these correspond to the usual braided
morphisms of crossed-modules (which, in turn, corresponds to the
braided morphisms of 2-groups under the equivalence of
$\S$\ref{SS:Equiv}).

It is easy to verify that two braided butterflies compose to a
braided butterfly. Braided crossed-modules and braided butterflies
form a bicategory $\mathcal{BRCM}$. There is a forgetful bifunctor
$\mathcal{BRCM} \to \mathcal{CM}$. The braided version of Theorem
\ref{T:bicat} is also true.

\subsection{Kernels and cokernels of braided butterflies}

The most interesting aspect of braided butterflies is that the NE-SW
of such a butterfly has a natural structure of a 2-crossed-module.

Also, the kernel, the cokernel, and the image  of a braided
butterfly $\mathcal{P} \: \bbH \to \bbG$ (see
$\S$\ref{SS:kernelbutterfly}) between braided crossed-modules are
naturally braided crossed-modules, and the  natural maps
$\Ker\mathcal{P} \to \bbH$ and $\bbG \to  \Coker\mathcal{P}$ are
braided morphisms of crossed-modules. The action of $G_2$ on $E$
given by  $x^g:=x\iota(\{\rho(x),g\}_{\bbG})$.

\subsection{Butterflies in an abelian category}{\label{SS:abelian}}

Let $\A$ be an abelian category. The results of the previous
sections are also valid for the category of complexes of length 2 in
$\A$. More precisely, let $\Ch(\A)$ be the bicategory whose objects
are complexes $X^{-1} \to X^0$, whose morphisms are
abelian\footnote{Abelian means that, since there are no actions, we
are dropping  the requirement for compatibility of actions.}
butterflies

 $$\xymatrix@C=8pt@R=6pt@M=6pt{ X^{-1} \ar[rd]^{\kappa} \ar[dd]_{\partial}
                            & & Y^{-1} \ar[ld]_{\iota} \ar[dd]^{\partial} \\
                  & E \ar[ld]^{\sigma} \ar[rd]_{\rho}  & \\
                  X^0 & & Y^0       }$$
and whose 2-morphisms are isomorphisms  of butterflies. The
composition of butterflies is defined as in the case of usual
butterflies ($\S$\ref{S:Bicat}). The criterion for strictness of a
butterfly is also valid ($\S$\ref{SS:st}).

Let us describe the additive structure on $\Ch(\A)$. Given two
morphisms $P$, $P'$ in $\Ch(\A)$
       $$\xymatrix@C=8pt@R=6pt@M=6pt{ X^{-1} \ar[rd]^{\kappa}
       \ar[dd]
                            & & Y^{-1} \ar[ld]_{\iota} \ar[dd]  \\
                  & E \ar[ld]^{\sigma} \ar[rd]_{\rho}  & \\
                  X^0 & & Y^0       } \ \   \ \
       \xymatrix@C=8pt@R=6pt@M=6pt{ X^{-1} \ar[rd]^{\kappa'} \ar[dd]
                            & & Y^{-1} \ar[ld]_{\iota'} \ar[dd] \\
                  & E' \ar[ld]^{\sigma'} \ar[rd]_{\rho'}  & \\
                  X^0 & & Y^0       }$$
with the same source and target, we define $P+P'$ to be the
butterfly
   $$\xymatrix@C=8pt@R=6pt@M=6pt{ X^{-1} \ar[rd]^{(\kappa,\kappa')}
   \ar[dd]
                            & & Y^{-1} \ar[ld]_{(0,\iota)} \ar[dd] \\
                  &
   E\oux{X^0}{Y^{-1}}E' \ar[ld]^{\sigma=\sigma'} \ar[rd] _{\rho+\rho'}  & \\
                  X^0 & & Y^0       }$$
We define $-P$ to be the butterfly
  $$\xymatrix@C=8pt@R=6pt@M=6pt{ X^{-1} \ar[rd]^{\kappa} \ar[dd]_{d}
                            & & Y^{-1} \ar[ld]_{-\iota} \ar[dd]^{d} \\
                  & E \ar[ld]^{\sigma} \ar[rd]_{-\rho}  & \\
                  X^0 & & Y^0       }$$

Using the above addition, the morphism groupoids in $\Ch(\A)$ become
symmetric monoidal categories, and the composition rule becomes a
symmetric monoidal functor.

The category obtained by identifying 2-isomorphic morphisms in
$\Ch(\A)$ is naturally equivalent to the full subcategory of the
derived category $\mathcal{D}(\A)$ consisting of complexes sitting
in degrees $[-1,0]$. Under this equivalence, the diagonal sequence
$$X^{-1}\llra{\kappa}E\llra{\rho}Y^0,$$
viewed as a complex sitting in degrees $[-2,0]$, corresponds to the
mapping cone.

 \begin{rem}
  As in the non-abelian case, an abelian butterfly as above comes from
  a chain map if and only if the NE-SW sequence is split. This is
  automatically the case if either $X_0$ is projective or $Y^{-1}$ is
  injective.
 \end{rem}

Let us define $\Chst(\A)$ to be the category whose objects are
complexes sitting in degrees $[-1,0]$, whose morphisms are morphisms
of complexes, and whose 2-morphisms are chain homotopies. Note that
morphism groupoids in $\Chst(\A)$ are (strict) symmetric monoidal.
There is natural (weak) functor $\Chst(\A) \to \Ch(\A)$ defined as
in $\S$\ref{SS:st}. (This is the 2-categorified version of the
quotient functor from the homotopy category to the derived
category.) This functor is a bijection on objects and faithful (and
symmetric monoidal) on morphism groupoids; see Theorem
\ref{T:bicat}.

Let $\mathbb{X}$ and $\mathbb{Y}$ be objects in $\Ch(\A)$. The next
proposition gives us some information about the groupoid
$\Hom_{\Ch(\A)}(\mathbb{X},\mathbb{Y})$.

 \begin{prop}{\label{P:exactseq}}
  Let $\mathbb{X}=[X^{-1}\to X^0]$ and $\mathbb{Y}=[Y^{-1}\to Y^0]$ be
  objects in $\Ch(\A)$. Let $\mathcal{E}xt(X^0,Y^{-1})$ be the
  groupoid whose objects are extensions of $X^0$ by $Y^{-1}$ and whose
  morphisms are isomorphisms of extensions. Then, the forgetful map
  $\Hom_{\Ch(\A)}(\mathbb{X},\mathbb{Y}) \to
  \mathcal{E}xt(X^0,Y^{-1})$ which sends a butterfly to its NE-SW
  sequence is a fibration of groupoids whose fiber is equivalent to
  $\Hom_{\Chst(\A)}(\mathbb{X},\mathbb{Y})$. Indeed, the sequence
    $$0 \to \Hom_{\Chst(\A)}(\mathbb{X},\mathbb{Y}) \to
     \Hom_{\Ch(\A)}(\mathbb{X},\mathbb{Y}) \to \mathcal{E}xt(X^0,Y^{-1})$$
  is an exact sequence of symmetric monoidal categories.
 \end{prop}

\section{A special case of Theorem \ref{T:butterfly}; Dedecker's theorem}
{\label{S:special}}

In this section we look at the special case of Theorem
\ref{T:butterfly} with $\bbH=[1\to\Gamma]$, where $\Gamma$ is a
group. The groupoid $\M(\Gamma,\bbG)$ looks as follows. The object
  of $\M(\Gamma,\bbG)$ are diagrams of the form
   $$\xymatrix@C=-4pt@R=10pt@M=10pt{
         & G_2\ar[ld]\ar[rd]^{\partial}    &   \\
       E \ar[rr]_(0.45){\rho}                  &     &   G_1         }$$
Here $E$ is an extension of $\Gamma$ by $G_2$, and for every $x \in
E$ and $\alpha \in G_2$ we require that
$\alpha^{\rho(x)}=x^{-1}\alpha x$. In fact, the short exact sequence
$$1 \to G_2 \to E \to \Gamma \to 1$$ is also part of the data, but
we suppress it from the notation and denote such a diagram simply by
$(E,\rho)$.

A morphism in $\M(\Gamma,\bbG)$ from $(E,\rho)$ to $(E',\rho')$ is
an isomorphism $f \: E\to E'$ of extensions (so it induces identity
on $G_2$ and $\Gamma$) such that $\rho=\rho'\circ f$.

 \begin{thm}{\label{T:classification}}
      The functor
          $$\Omega \:\Map(\Gamma,\bbG) \to  \M(\Gamma,\bbG)$$
      is an equivalence of groupoids. That is, $\M(\Gamma,\bbG)$
      is naturally equivalent to the groupoid of weak functors from
      $\Gamma$ to $\bbG$.
 \end{thm}

 \begin{proof}
   Follows immediately from Theorem \ref{T:butterfly}.
 \end{proof}

 \begin{cor}[Dedecker, \oldcite{Ded}]{\label{C:homotopy}}
    There is a natural bijection
      $$\pi_0 \: \hhom(\Gamma,\bbG) \risom \pi_0 \M(\Gamma,\bbG).$$
    In other words, the   homotopy
    classes of weak maps from $\Gamma$ to $\bbG$ are in a natural
    bijection with isomorphism classes of diagrams of the form
         $$\xymatrix@C=-4pt@R=10pt@M=10pt{
            & G_2\ar[ld]\ar[rd]^{\partial}    &   \\
             E \ar[rr]_(0.45){\rho}                  &     &   G_1         }$$
     where the $E$ is an extension of $\Gamma$ by $G_2$, and for every $x \in E$
     and $\alpha \in G_2$ the equality $\alpha^{\rho(x)}=x^{-1}\alpha x$ is satisfied.
 \end{cor}


 \begin{ex}{\label{E:extension}}
  Let $A$ be an abelian group, and consider the crossed-module
  $\mathfrak{A}=[A \to 1]$. The groupoid of central extensions of
  $\Gamma$ by $A$ is naturally equivalent to the groupoid
  $\M(\Gamma,\mathfrak{A})$, which is itself naturally equivalent to
  the derived mapping groupoid $\Map(\Gamma,\mathfrak{A})$, by
  Theorem \ref{T:classification}. In particular, by Corollary
  \ref{C:homotopy}, the set of homotopy classes of weak maps from
  $\Gamma$ to  $\mathfrak{A}$ is in natural bijection with
  isomorphism  classes of central extensions of $\Gamma$ by $A$,
  which is itself in natural bijection with
  $H^2(\Gamma,A)$.\footnote{This is of course not surprising, since
  $\mathcal{A}$ is an algebraic model for the Eilenberg-MacLane
  space $K(A,2)$.}
 \end{ex}

\begin{ex}{\label{E:aut}}
  Let $K$ be a group. Let $\mathfrak{Aut}(K)$ be the crossed-module
  $[K \to \Aut(K)]$.\footnote{The corresponding 2-group is the
  2-group of self-equivalences of $K$, where $K$ is viewed as a
  category with one object.} The groupoid of extensions of $\Gamma$
  by $K$ is  naturally equivalent to the groupoid
  $\M(\Gamma,\mathfrak{Aut}(K))$, which is itself naturally
  equivalent to the derived mapping groupoid
  $\Map(\Gamma,\mathfrak{Aut}(K))$, by Theorem
  \ref{T:classification}. In particular, by Corollary
  \ref{C:homotopy}, the set of homotopy classes of weak maps from
  $\Gamma$ to  $\mathfrak{Aut}(K)$ is in natural
  bijection with isomorphism  classes of extensions of $\Gamma$ by $K$.
 \end{ex}

 \begin{rem}
  The above example is identical to Theorem 2 (and Proposition 3) of
  \cite{BBF} in the case where $\mathcal{G}$ (notation as in {\em
  loc.\;cit.}) has only one object. It seems that, indeed, the general
  case of Theorem 2 (and also Proposition 3) of \cite{BBF}, with
  $\mathfrak{G}$ an arbitrary groupoid,  is equivalent to this special
  case. To see this, note that we may assume $\mathcal{G}$ is
  connected. On the other hand, since both sides of the equality in
  Theorem 2 (and also Proposition 3) of \cite{BBF} are functorial in
  $\mathcal{G}$, we may assume $\mathcal{G}$ has only one object.
 \end{rem}

 \begin{rem}
  Corollary \ref{C:homotopy} was first discovered by Dedecker
  \oldcite{Ded}. It is discussed in detail in \cite{Breen} in the more
  general case where everything takes place in a Grothendieck site. In
  analogy with Example \ref{E:aut}, an element in $\pi_0
  \M(\Gamma,\bbG)$ is called a {\em $\bbG$-extension of $\Gamma$}.
  Such extensions can be thought of as ``$\bbG$-torsors'' over the
  classifying space $B\Gamma$,  the same way that  group extensions of
  $\Gamma$ by $K$ can be regarded as $\mathfrak{Aut}(K)$-torsors
  (i.e., $K$-gerbes) over $B\Gamma$. Therefore, one can alternatively
  describe elements of $\pi_0 \M(\Gamma,\bbG)$ as certain first
  cohomology classes of $\Gamma$ with coefficients in the
  crossed-module $\bbG$. For more on this, the reader is encouraged to
  consult \cite{Breen}, especially, $\S$8.
 \end{rem}

In the rest of this section we present three isolated facts for
which we have no use, but we found them interesting nevertheless!

\subsection{Side note 1: split crossed-modules}
We present a cute (and presumably well-known) application of
Corollary \ref{C:homotopy}.

A 2-group $\mathfrak{G}$ is called {\em split} if it is completely
determined by  $\pi_1\mathfrak{G}$, $\pi_2\mathfrak{G}$, and the
action of the former on the latter. More precisely, $\mathfrak{G}$
is split if it is isomorphic in $\Ho(\cm)$ to the crossed-module
$[\partial \: \pi_2\mathfrak{G} \to \pi_1\mathfrak{G}]$, where
$\partial$ is the trivial homomorphism. From the homotopical point
of view, the following proposition is straightforward (see the
beginning of $\S$\ref{S:Homotopical}). However, to give a purely
algebraic  proof seems to be tricky. We give a proof that makes use
of Corollary \ref{C:homotopy}.

  \begin{prop}{\label{P:split}}
     Let $\mathfrak{G}=[G_2 \to G_1]$ be a crossed-module, and
     assume that the map $\bbG \to \pi_1\bbG$, viewed in
     $\Ho(\cm)$, admits a section.
     Then $\mathfrak{G}$ is split.
  \end{prop}

  \begin{proof}
   By Corollary \ref{C:homotopy}, there exists an extension $$
   1\arr{} G_2 \arr{} E \arr{f} \pi_1\bbG\arr{} 1$$ and a map $\rho
   \: E \to G_1$ satisfying the conditions stated therein. Consider
   the semi-direct product $G_2\ltimes\pi_2\bbG$ where $G_2$ acts on
   $\pi_2\bbG$ by conjugation. It fits in a crossed-module
   $\bbG'=[G_2\ltimes\pi_2\bbG \to E]$ where the map
   $G_2\times\pi_2\bbG \to E$ is obtained from the first projection
   map, and the action of $E$ on both factors is by conjugation. We
   have a homomorphism $\sigma \: G_2\ltimes\pi_2\bbG \to G_2$ which
   on the first factor is just the identity and on the second factor
   is the inclusion map. It is easy to see that $(\sigma,\rho) \:
   \bbG' \to \bbG$ is a crossed-module map; in fact, it is an
   equivalence of crossed-modules. On the other hand, $(pr_2,f) \:
   [G_2\ltimes\pi_2\bbG \to E] \to [\pi_2\bbG \to \pi_1\bbG]$ is also
   an equivalence of crossed-modules. So we have constructed a zigzag
   of equivalences that connect $\bbG$ to the trivial crossed-module
   $[\pi_2\bbG \to \pi_1\bbG]$. So $\bbG$ is split.
 \end{proof}

\subsection{Side note 2: Hopf's formula}
It is interesting  to compute $\pi_0\Map(\Gamma,\mathfrak{A})$
straight from  definition of the derived mapping groupoid
(Definition \ref{D:mapping}). This leads to Hopf's formula for
$H^2(\Gamma,A)$.

Choose a presentation $F/R\cong\Gamma$ of $\Gamma$, where $F$ is
free. The crossed-module maps  $[R\to F] \to [A\to 1]$ are precisely
the group homomorphism $g \: R \to A$ that are constant on the
conjugacy classes (under the $F$-action) of elements in $R$. In
other words, $g(x)=1$ for every $x \in [F,R]$. Two such
homomorphisms $g$ and $g'$ are homotopic, if there is a group
homomorphism $h \: F \to A$ such that $h|_R=g'g^{-1}$.  So,
$\pi_0\Map(\Gamma,\mathfrak{A})$ is in natural bijection with
$$\Coker\left\{\Hom(F,A) \to \Hom(R/[F,R],A)\right\}.$$ This is
Hopf's famous formula for $H^2(\Gamma,A)$.

\section{Cohomological point of view}{\label{S:Cohomological}}

In this section we give a cohomological characterization of
$\mathcal{M}(\Gamma,\bbG)$, Theorem \ref{T:maps}, and compare it to
Theorem \ref{T:classification}. Theorem \ref{T:maps} is
well-known\footnote{The diligent reader can dig it out of
\oldcite{AzCe}, $\S$5.}, but the explicit way in which it relates to
Theorem \ref{T:classification} is what we are interested in. Since
this construction has been quoted in \cite{B-N}, we feel obliged to
include the precise account. We begin by recalling some standard
facts about group extensions.

\subsection{Groups extensions and $H^2$; review}{\label{SS:Extensions}}

We recall  some basic facts about classification of group
extensions via cohomological invariants. Our main reference is
\cite{Brown}.\footnote{Also see Example \ref{E:aut}.}

Let $N$ and $\Gamma$ be groups (not necessarily abelian). We would
like to classify extensions
  $$1 \to N \to E \to \Gamma \to 1,$$
up to isomorphism. First of all, notice that such an extension gives
rise to a group homomorphism $\psi \: \Gamma \to \Out(N)$. So we
might as well fix $\psi$ as part of the data, and classify extension
which induce $\psi$. Denote the set of such extensions by
$\mathcal{E}(\Gamma,N,\psi)$. Let $C$ denote the center of $N$, made
into a $\Gamma$-module through $\psi$. We have the following
theorem.

 \begin{thm}[\oldcite{Brown}, Theorem 6.6]\label{T:extensions}
  The set $\mathcal{E}(\Gamma,N,\psi)$ admits a natural free,
  transitive action by the abelian group $H^2(\Gamma,C)$. Hence,
  either $\mathcal{E}(\Gamma,N,\psi)=\emptyset$, or else there is a
  bijection $\mathcal{E}(\Gamma,N,\psi)\leftrightarrow
  H^2(\Gamma,C)$. This bijection
  depends on the choice of a particular element of
  $\mathcal{E}(\Gamma,N,\psi)$.
 \end{thm}

 \begin{rem}{\label{R:extensions}}
  If we are given a lift $\tilde{\psi} \: \Gamma \to \Aut(N)$ of
  $\psi$,  we obtain a distinguished element in
  $\mathcal{E}(\Gamma,N,\psi)$, namely, the semi-direct product
  $N\rtimes \Gamma$. This is automatically the case if, for
  instance, $N$ is abelian. Therefore, when $N$ is abelian, we have
  a canonical bijection
  $\mathcal{E}(\Gamma,N,\psi)\leftrightarrow H^2(\Gamma,C)$.
 \end{rem}

The meaning of this theorem is that, given two elements $E_0$, $E$
in $\mathcal{E}(\Gamma,N,\psi)$, one can produce their {\em
difference} as an element in $H^2(\Gamma,C)$. Notice that $C$ is now
abelian, so, by Remark \ref{R:extensions}, every element in
$H^2(\Gamma,C)$ gives rise to a canonical extension of $\Gamma$ by
$C$. Below we will explain how this extension can be explicitly
constructed from $E_0$ and $E$.

Let us call a complex $ M \to E \to \Gamma$ {\em semi-exact} if the
left map is injective, the right map is surjective, and  the kernel
$K$ of  $E \to \Gamma$ is generated by  $M$ and $C_K(M)$ (the the
centralizer of $M$ in $K$). This last condition guarantees that
there is a well-defined homomorphism $\psi \: \Gamma \to \Out(M)$.

Assume we are given two semi-exact sequences
      $$ M \to E_0 \to \Gamma  \ \ , \ \ M \to E \to \Gamma$$
such that $\psi_0=\psi$. Define $L$ to be the group of pairs $(x,y)
\in E_0\times E$  such that $\bar{x}=\bar{y} \in \Gamma$, and that
conjugation by $x$ and $y$ induce the same automorphism of $M$.
Observe that $I=\{(a,a) \in E_0\times E\ | \ a\in M\}$ is a normal
subgroup of $L$.

 \begin{defn}[Baer product]{\label{D:product}}
   Notation being as above, define $E_0\oux{\Gamma}{M}E:=L/I$.
 \end{defn}

There is an obvious surjective homomorphism $E_0\oux{\Gamma}{M}E \to
\Gamma$. It fits in a natural  exact sequence
    $$1 \to C \to C_{K_0}(M)\times C_{K}(M) \to E_0\oux{\Gamma}{M}E
               \to \Gamma \to 1,$$
where $C$ is the center of $M$ mapping diagonally to  $C_K(M)\times
C_{K'}(M)$.  (This gives us two semi-exact sequences:
    $$C_{K_0}(M) \to E_0\oux{\Gamma}{M}E
               \to \Gamma \ \ , \ \ C_{K}(M) \to E_0\oux{\Gamma}{M}E
               \to \Gamma$$
which are somehow mirror to each other.)

If the two semi-exact sequences that we started with were actually
exact, we would have $C_K(M)=C_{K'}(M)=C$. So, by identifying $C$ as
the cokernel of the diagonal map $C \to C\times C$, we obtain the
following exact sequence $$1 \to C  \to E_0\oux{\Gamma}{M}E      \to
\Gamma \to 1.$$ More explicitly, we define the map $C \to
E_0\oux{\Gamma}{M}E$  by sending  $a$ to $(a,1)$.

 \begin{defn}{\label{D:difference}}
  Let $E_0,E \in \mathcal{E}(\Gamma,N,\psi)$.
  Define the {\em difference}  $D(E_0,E)$ to be the
  the sequence
  $1 \to  C \to E_0\oux{\Gamma}{N}E \to \Gamma\to 1$
  defined above.
 \end{defn}

 \begin{rem}{\label{R:symmetric}}
  Observe that $E_0\oux{\Gamma}{M}E$ is symmetric with respect to
  $E_0$ and $E$, but $D(E_0,E)$ is not. What determines the sign in
  the above construction is the map $C \to E_0\oux{\Gamma}{N}E$. So,
  if instead of  $a \mapsto (a,1)$ we used $a \mapsto (1,a)$ we
  would obtain $D(E,E_0)$, because in $C \to E_0\oux{\Gamma}{M}E$
  the elements  $(1,a)$ and $(a,1)$ are inverse to each other.
 \end{rem}

Conversely, given  $E_0 \in  \mathcal{E}(\Gamma,N,\psi)$ and an
extension   $ 1 \to C \to H \to \Gamma \to 1$ (recall that $C$ is
the center of $N$), we can recover $E$ as the difference
$E:=D(E_0,H)$. In other words, consider the group
$E_0\oux{\Gamma}{C}H$. This contains $N=N\oux{1}{C}C$ as a normal
subgroup. The sequence
$$  1 \to N \to E_0\oux{\Gamma}{C}H \to \Gamma \to 1$$
is the desired extension.

\subsection{Cohomological classification of maps into a 2-group}
{\label{SS:CohClass}}

In this subsection we prove the following cohomological
classification of the homotopy classes of weak maps from a group
$\Gamma$ to a 2-group $\bbG$. The result itself is well-known, but
the way in which it relates to the classification theorem (Theorem
\ref{T:classification}) is what we are interested in.

  \begin{thm}{\label{T:maps}}
     Let $\mathfrak{G}$ be a 2-group, and let $\Gamma$ be a discrete
     group. Fix a homomorphism $\chi \: \Gamma \to
     \pi_1\mathfrak{G}$, and let
     $[\Gamma,\mathfrak{G}]_{\tgp}^{\chi}$ be the set of
     homotopy classes of weak maps $\Gamma \to \bbG$ inducing $\chi$
     on $\pi_1$. Then, either
     $[\Gamma,\mathfrak{G}]_{\tgp}^{\chi}$ is empty, or  it
     is naturally a transitive $H^2(\Gamma,\pi_2\mathfrak{G})$-set.
     (Here, $\pi_2\mathfrak{G}$ is made into a $\Gamma$ module via
     $\chi$.)
  \end{thm}

In $\S$\ref{S:Homotopical} we see exactly when
  $[\Gamma,\mathfrak{G}]_{\tgp}^{\chi}$
  is non-empty (see Remark \ref{R:vanish}) .

\begin{rem}
  If in the above proposition we take $\bbG=\mathfrak{Aut}(N)$
  we recover Theorem \ref{T:extensions}; see Example \ref{E:aut}.
\end{rem}

Before proving  Theorem \ref{T:maps} we need some preliminaries.

\vspace{0.1in} \noindent{\em Conventions for this subsection.}
Throughout this subsection, we fix $\chi$ (and consequently an
action of $\Gamma$ on $\pi_2\bbG$). We will think of
$H^2(\Gamma,\pi_2\mathfrak{G})$ as the group of isomorphism classes
of extensions of $\Gamma$ by $\pi_2\mathfrak{G}$ for which the
induced action of $\Gamma$ on $\pi_2\bbG$ is the one we have fixed.
Whenever we talk about an extension of $\Gamma$ by
$\pi_2\mathfrak{G}$ we assume that this condition is satisfied.
\vspace{0.1in}

\noindent{\bf Construction of the  action of
$H^2(\Gamma,\pi_2\mathfrak{G})$ on
$[\Gamma,\mathfrak{G}]_{\tgp}^{\chi}$.}

Suppose we are given $(E_0,\rho_0)$ in
$[\Gamma,\mathfrak{G}]_{\tgp}^{\chi}$ and an extension
   $$1 \to \pi_2\bbG \to K \to \Gamma \to 1$$
in $H^2(\Gamma,\pi_2\mathfrak{G})$. Set
$E:=E_0\underset{\Gamma}{\overset{\pi_2\mathfrak{G}}{\times}}K$ (see
Definition \ref{D:product}). The inclusion  $G_2 \hookrightarrow
E_0$ induces a natural homomorphism $G_2 \to E$ which identifies
$G_2$ with a normal subgroup of $E$. The quotient is $\Gamma$. We
have a natural map $\rho \: E \to G_1$ defined by
$\rho(x,a)=\rho_0(x)$. This is easily seen to be well-defined.
Finally, it is easy to check that the   action of $E$ on $G_2$
induced via $\rho$ is equal to the conjugation action of $E$ on
$G_2$. This gives us the desired diagram:
   $$\xymatrix@C=-4pt@R=10pt@M=10pt{
         & G_2\ar[ld]\ar[rd]^{\partial}    &   \\
       E \ar[rr]_(0.45){\rho}        &     &   G_1         }$$

\vspace{0.1in}

\noindent{\bf Transitivity of the action.} To prove the transitivity
of the action of $H^2(\Gamma,\pi_2\mathfrak{G})$ on
$[\Gamma,\mathfrak{G}]_{\tgp}^{\chi}$ we employ a `difference
construction' similar to the ones of the previous subsection. It
takes  two elements in $[\Gamma,\mathfrak{G}]_{\tgp}^{\chi}$ and
produces their difference, which is  an element in
$H^2(\Gamma,\pi_2\mathfrak{G})$.

Assume $[\Gamma,\mathfrak{G}]_{\tgp}^{\chi}$ is non-empty, and fix
an   element $(E_0,\rho_0)$ in it as in the following diagram (see
Theorem \ref{T:classification}):
   $$\xymatrix@C=-4pt@R=10pt@M=10pt{
         & G_2\ar[ld]\ar[rd]^{\partial}    &   \\
       E_0 \ar[rr]_(0.45){\rho_0}                  &     &   G_1         }$$
Let  $(E,\rho)$ be another such diagram. Define the group $L$ by
   $$L=\{(x,y) \in E_0\times  E \ | \
               \bar{x}=\bar{y},\, \rho_0(x)=\rho(y) \}.$$
There is a natural surjective group homomorphism $L \to \Gamma$
sending $(x,y)$ to $\bar{x}=\bar{y}$. The kernel of this map is the
following group:
 $$\{(\alpha,\beta) \in G_2\times G_2 \
                | \ \alpha\beta^{-1}  \in
                      \pi_2\mathfrak{G}\}=I\times \pi_2\mathfrak{G}$$
where
   $$I:=\{(\beta,\beta) \in E_0\times E \
                | \ \beta \in G_2\},$$
and $\pi_2\bbG$ is identified with the subgroup of elements of the
form $(\alpha,1)$, $\alpha \in \pi_2\bbG$. It is  easy   to check
that $I$ is normal in $L$.

 \begin{defn}{\label{D:difference2}}
  Define $E_0\times_{\bbG}E=L/I$. The map $\alpha \mapsto
  (\alpha,1)$ identifies $\pi_2\bbG$ with a normal subgroup of $L$
  with cokernel $\Gamma$. The   extension $$1 \to \pi_2\bbG \to
  E_0\times_{\bbG}E \to \Gamma \to 1,$$ or its class in
  $H^2(\Gamma,\pi_2\bbG)$, is called the {\em difference} of
  $(E_0,\rho_0)$ and $(E,\rho)$ and is denoted by
  $D\big((E_0,\rho_0),(E,\rho)\big)$.\footnote{ Definition
  \ref{D:difference} is a special case of Definition
  \ref{D:difference2} with $\bbG=\mathfrak{Aut}(N)$; see
  Example \ref{E:aut}.}
 \end{defn}

 \begin{proof}[Proof of Theorem \ref{T:maps}]
  We leave it to the reader to verify that the construction of the
  action of $H^2(\Gamma,\pi_2\mathfrak{G})$ on
  $[\Gamma,\mathfrak{G}]_{\tgp}^{\chi}$
  and the difference construction of Definition \ref{D:difference2}
  are inverse to each other.
 \end{proof}

An interesting special case is when $\chi \: \Gamma \to \pi_1\bbG$
can be lifted to $\tilde{\chi} \: \Gamma \to G_1$. In the following
corollary we fix such a lift.

 \begin{cor}{\label{C:lifts}}
  With the hypothesis of the preceding paragraph, every class $f \in
  [\Gamma,\mathfrak{G}]_{\tgp}^{\chi}$ is uniquely
  characterized by (the isomorphism class of)  an extension $$1 \to
  \pi_2\mathfrak{G} \to K \to \Gamma  \to 1.$$ More explicitly,
  given such an extension we obtain $(E,\rho)$, where $E:=K
  \ltimes^{\pi_2\mathfrak{G}} G_2$ (Definition \ref{D:crossed}), and
  $\rho(k,\alpha):=\tilde{\chi}(\bar{k})\underline{\alpha}$, for
  $(k,\alpha) \in E$. Here the action of $K$ on $G_2$
  is obtained via $\tilde{\chi}$ from that of $G_1$ on $G_2$.
 \end{cor}

 \begin{proof}
   Set $E_0=\Gamma\ltimes G_2$, the action being obtained through
  $\tilde{\chi}$, and define $\rho_0 \: E_0 \to G_1$ by
  $\rho(x,\alpha):=\chi(x)\underline{\alpha}$. (Keep in kind  that
  $K\underset{\Gamma}
  {\overset{\pi_2\mathfrak{G}}{\times}}(\Gamma\ltimes G_2)=
  K \ltimes^ {\pi_2\mathfrak{G}}G_2$.)
 \end{proof}

An important special case of the above corollary is when there
exists a section for the map $G_1 \to \pi_1\mathfrak{G}$. In this
case, after fixing such a section, we have an explicit description
of the elements in $[\Gamma,\mathfrak{G}]_{\tgp}$, for any group
$\Gamma$.

\section{Homotopical point of view}{\label{S:Homotopical}}

The results of the previous sections are best understood if viewed
from the point of  homotopy theory of 2-types. Recall that 2-groups
(respectively, weak functors between 2-groups, weak natural
transformations) are algebraic models for pointed and connected
homotopy 2-types (respectively, pointed continuous maps, pointed
homotopies). Using this dictionary, the problem of classification of
weak maps between 2-groups translates to the problem of
classification of (pointed) homotopy classes of continuous maps
between (pointed) homotopy 2-types. The latter can be solved using
standard technique from obstruction theory.

In this section we explain how the homotopical approach works and
give a homotopical proof of Theorem \ref{T:maps}. In the next
section we compare the homotopical approach with the approaches of
the previous section. These all are presumably folklore and
well-known.

\subsection{The classifying space functor and  homotopical proof of
Theorem \ref{T:maps}}{\label{SS:Homotopical}}

Consider the {\em classifying space} functor $B \: \Ho(\tgp) \to
\Ho(\mathbf{Top}_*)$ defined by $B(\bbG)=|N(\bbG)|$. By Corollary
\ref{C:equiv} of Appendix, we know that this gives rise to an
equivalence between the homotopy category of 2-groups and the
homotopy category of pointed connected CW-complexes with vanishing
$\pi_i$, $i\geq 3$ (the so called {\em connected homotopy 2-types}),
a result essentially due to Whitehead \cite{W}. So it would be most
natural to   view the algebraic results of the previous sections
from this homotopic perspective. To complete the picture, we  recall
MacLane-Whitehead  characterization of pointed connected homotopy
2-types from \cite{M-W}.

Mac Lane and Whitehead show that, to give a connected pointed
homotopy 2-types is equivalent to giving a triple
$(\pi_1,\pi_2,\kappa)$ where $\pi_1$ is an arbitrary group, $\pi_2$
is an abelian group endowed with a $\pi_1$ action, and $\kappa \in
H^3(\pi_1,\pi_2)$. To see where such a triple comes from, let $X$ be
a pointed connected CW-complex such that $\pi_iX=0$, $i\geq 3$.
Consider the Postnikov decomposition of $X$:
  $$\xymatrix@R=12pt@C=-5pt@M=6pt{
      K(\pi_2,2) \ - \ar@{-}  & X_2  \ar[d] & \\
                      &  X_1& =K(\pi_1,1)     }$$
The triple corresponding to $X$ is $(\pi_1X,\pi_2X,\kappa)$, where
$\kappa$ is the Postnikov invariant corresponding to the above
picture. Recall that this  Postnikov invariant is the obstruction to
existence of a section for the fibration $p \: X_2 \to X_1$. We know
from obstruction theory that, if this obstruction vanishes, then for
any choice of base points $x_1 \in X_1$ and $x_2 \in p^{-1}(x_1)$,
there exist a  pointed section $X_1 \to X_2$. Furthermore, after
fixing such a section, the pointed homotopy classes\footnote{We
should actually be considering {\em fiberwise} homotopy classes, but
since in our case the fiber is simply connected we get the same
thing.} of such section are in bijection with $H^2(\pi_1X,\pi_2X)$.

 \begin{proof}[Homotopical proof of Theorem \ref{T:maps}]
    Passing to classifying space induces a bijection between
    $[\Gamma,\mathfrak{G}]_{\tgp}$
    and  pointed homotopy class of
    maps $B\Gamma \to B\mathfrak{G}$. Let
        $$\xymatrix@R=12pt@C=-5pt@M=6pt{
      K(\pi_2\bbG,2) \ - \ar@{-}  & X_2  \ar[d] & \\
                      &  X_1& =K(\pi_1\bbG,1)     }$$
    be the Postnikov tower of $B\mathfrak{G}$. To give a group
    homomorphism $\chi \: \Gamma \to \pi_1\bbG$ is equivalent to
    giving a pointed homotopy class from $B\Gamma$ to $X_1$. Fix
    such a class $F \: B\Gamma \to X_1$. The question is now to
    classify the pointed (equivalently, fiberwise -- because the
    fiber is simply connected) homotopy classes of lifts $f \:
    B\Gamma \to X_2$ of $F$. We know from obstruction theory that
    the obstruction to existence of such a lift is precisely
    $F^*(\kappa)=\chi^*(\kappa)\in H^3(\Gamma,\pi_2X)$, where
    $\pi_2X$ is made into a $\Gamma$-module via $\chi$. By
    obstruction theory, whenever this obstruction vanishes, the set
    of pointed (equivalently, fiberwise) homotopy classes of lifts
    of such lifts $f$ is a transitive     $H^2(\Gamma,\pi_2X)$-set.
  \end{proof}

 \begin{rem}{\label{R:vanish}}
  As we saw in the above proof, for a fixed $\chi \: \Gamma \to
  \pi_1\bbG$, the set $[\Gamma,\mathfrak{G}]_{\tgp}^{\chi}$
  is non-empty if and only if $\chi^*(\kappa)\in
  H^3(\Gamma,\pi_2\bbG)$ is zero, where $\kappa$ is the Postnikov
  invariant  of $\bbG$.
 \end{rem}

 \begin{ex}{\label{E:aut2}}
  Let $\Gamma$ and $N$ be discrete groups. For a given $\chi \:
  \Gamma \to \Out(N)$, the obstruction to lifting $\chi$ to a map
  $\Gamma \to  \mathfrak{Aut}(N)$ (equivalently, to finding an
  extension of $\Gamma$ by $N$ giving giving rise to $\chi$ -- see
  Example \ref{E:aut}), is the element $\chi^*(\kappa) \in
  H^3(\Gamma,C(N))$, where $C(N)$ is the center of $N$ and $\kappa
  \in H^3(\Out(N),C(N))$  is the Postnikov invariant of
  $\mathfrak{Aut}(N)$. (In other words, the Postnikov invariant
  $\kappa$ of $\mathfrak{Aut}(N)$ is the universal obstruction class
  for the existence of group extensions.) When $\chi^*(\kappa)=0$,
  the set of lifts of $\chi$ to $\mathfrak{Aut}(N)$ (equivalently,
  extensions of $\Gamma$ by $N$ giving rise to $\chi$)  admits a
  natural transitive action of $H^2(\Gamma,C)$, where
  $C=\pi_2\mathfrak{Aut}(N)$ is the  center of $N$.
 \end{ex}

For the sake of amusement, we also include a homotopical proof for
Proposition \ref{P:split}.

 \begin{proof}[Homotopical proof of Proposition \ref{P:split}]
  Let $\bbG$ be as in Proposition \ref{P:split}. We have to show
  that the 2-type $B\bbG$ is split. That is, it is homotopy
  equivalent to the product $K(\pi_1\mathfrak{G},1)\times
  K(\pi_2\mathfrak{G},2)$ of Eilenberg-MacLane spaces. We know that
  $\mathfrak{G}$ is classified by the triple
  $(\pi_1\mathfrak{G},\pi_2\mathfrak{G},\kappa)$. By assumption, the
  action of $\pi_1\mathfrak{G}$ on $\pi_2\mathfrak{G}$ is trivial.
  Also, by assumption, the corresponding Postnikov tower $X_2 \to
  X_1$ has a section, so $\kappa$  vanishes. Since the 2-type
  $K(\pi_1\mathfrak{G},1)\times K(\pi_2\mathfrak{G},2)$ also gives
  rise to the same triple, it must be (pointed) homotopy
  equivalent to $B\bbG$, which is what we wanted to prove.
 \end{proof}

\section{Compatibility of different approaches}{\label{S:Compatible}}

Let $p\: Y \to X$ be a fibration of CW-complexes (or simplicial
sets). In this section we recall  the notion of {\em difference
fibration} for two liftings of a map $F \: A \to X$,  and use it to
clarify Definition  \ref{D:difference2}, as well as to explain why
the cohomological classification of maps into a 2-group
($\S$\ref{SS:CohClass}) is compatible with the homotopical approach
($\S$\ref{S:Homotopical}).

A more detailed discussion of the difference fibration
construction, and its application to obstruction theory, can be
found in \cite{Baues}.

\subsection{Difference fibration construction}{\label{SS:Difference1}}

The difference constructions of $\S$\ref{SS:Extensions} and
$\S$\ref{SS:CohClass} are special cases of (and were originally
motivated by) a general difference construction for maps of
simplicial sets (or topological spaces). In this section we review
this general construction and explain how it relates to the
algebraic versions of it that we have already encountered in
previous sections.

Let $X$ and $Y$ be simplicial sets and $p \: Y \to X$  a simplicial
map. Let $A$ be another simplicial set and $F \: A \to X$ a map of
simplicial sets.  We are interested in the classification of lifts
of $f$ to $Y$, if such lifts exist. In our case of interest,  $A$
and $X$ are going to be 1-types and $Y$ a 2-type, so everything is
explicit and easy.

A useful tool in the study of such lifting problems is the {\em
difference construction}, as in (\cite{Baues}, page
293).\footnote{We will use the simplicial version of Baues's
definition though.}

Let $f_0,f \: A \to Y$ be liftings of $F$. To measure the
difference between $f_0$ and $f$, we construct the simplicial set
$D_p(f_0,f)$ as in the following cartesian diagram:
  $$\xymatrix@=12pt@M=10pt{
      D_p(f_0,f) \ar[r]\ar[d]  &   \ar[d]^{(d_0,d_1,p^{\Delta^1})}
                                                  Y^{\Delta^1}\\
          A  \ar[r]_(0.35){(f_0,f,c)}        &
                                Y\times Y\times X^{\Delta^1}   }$$
Here,  $c \: A \to X^{\Delta^1}$ is the map that sends $a \in A$ to
the constant path at $F(a)$. We are actually interested in the case
where $p \: Y \to X$
is a fibration.

The following lemma justifies the terminology, but since we will
not need it here we will not give the proof (except  in a special
case -- see Proposition \ref{P:compatible}). In the case of
topological spaces, a proof is can be found in \cite{Baues}.

 \begin{lem}{\label{L:difference}}
   If $p  \: Y \to X$ is a fibration, then $D_p(f_0,f) \to A$ is also
   a fibration.
 \end{lem}

The usefulness of the difference construction is justified by the
following proposition.

 \begin{prop}{\label{P:section}}\noindent\par
   \begin{itemize}

    \item[$\mathbf{i.}$] The simplicial set of sections to the
      map $D_p(f_0,f) \to A$ is naturally isomorphic to the
      simplicial set of fiberwise homotopies  between $f_0$ and $f$.

    \item[$\mathbf{ii.}$]  The primary difference of the sections
     $f_0$ and $f$ is precisely the primary obstruction to the
     existence of a section to $D_p(f_0,f) \to A$.
   \end{itemize}
 \end{prop}

 \begin{proof}
  Part ($\mathbf{i}$) follows from the definition. Part
  ($\mathbf{ii}$) is proved in \cite{Baues} (see page 295 {\em
  loc.\;cit.}) in the case of topological
  spaces. The proof can be adopted to the simplicial situation.
 \end{proof}

When all spaces and  maps are pointed,  $D_p(f_0,f)$ is also
naturally pointed. In this case we have:

 \begin{prop}{\label{P:pointedsection}}\noindent\par
   \begin{itemize}
    \item[$\mathbf{i.}$] The simplicial set of pointed sections to the
      map $D_p(f_0,f) \to A$ is naturally isomorphic to the
      simplicial set of fiberwise homotopies between $f_0$ and $f$
      which fix the base point.

    \item[$\mathbf{ii.}$]  The primary difference of the sections
      $f_0$ and $f$ is precisely the primary obstruction to the
      existence of a section to $D_p(f_0,f) \to A$ (everything
      pointed).
   \end{itemize}
 \end{prop}

Difference fibration construction can indeed be performed in any
category with fiber products in which there is an interval, and we
have a notion of internal hom.\footnote{We are not asking for a
monoidal structure here.} For instance, $\mathbf{Cat}$,
$\mathbf{2Cat}$, and $\mathbf{2Gpd}$ are examples of such a
category, where for the interval we take $\mathbf{I}_1=\{0 \to
1\}$, and the internal homs are given by $\homc$, as in Definition
\ref{D:homo2}. So, given a diagram
   $$\xymatrix@=30pt@M=6pt{         &  \mfD \ar[d]^p  \\
            \mfA  \ar@/^/[ru]|-{f_0} \ar@/_/[ru]|-{f}  \ar[r]_F  &
                                                  \mfC    }$$
of 2-categories, we can talk about the difference $D_p(f_0,f)$ of
$f_0$ and $f$. This is a 2-category with a natural functor
$D_p(f_0,f) \to A$. Furthermore, if $\mfA$, $\mfC$ and $\mfD$ are
pointed (i.e. have a chosen object), then so is $D_p(f_0,f)$. The
following 2-categorical version of Proposition \ref{P:section} is
also valid.

 \begin{prop}{\label{P:2catsection}}
  Suppose $\mfA$, $\mfC$ and $\mfD$ are 2-categories (respectively,
  pointed 2-categories) as above. Then the 2-category of sections
  (respectively, pointed sections) to the functor $D_p(f_0,f) \to A$
  is naturally isomorphic to the 2-category of fiberwise
  transformations
  (respectively, pointed fiberwise transformations) between $f_0$ and $f$.
 \end{prop}

Clearly if $\mfA$,  $\mfC$ and $\mfD$ are 2-groupoids, then so is
$D_p(f_0,f)$. So, we can talk about difference construction for
2-groupoids. It is also true that, if $p$ is  a fibration (Appendix,
Definition \ref{D:2fibration}), then so is $D_p(f_0,f) \to A$. The
latter statement, whose simplicial counterpart  we did not prove
(Lemma \ref{L:difference}), can be proved easily by verifying the
conditions of Definition \ref{D:2fibration}.

 \begin{defn}{\label{D:difference2gp}}
  In the above situation, assume  $\mfA$, $\mfC$ and $\mfD$ are
  2-groupoids with one object (i.e. 2-groups), and let  $* \in \Ob
  D_p(f_0,f)$ be the canonical base point of $D_p(f_0,f)$.
  We define $D_p(f_0,f)_*$  to be the 2-group of automorphisms of the object $*$.
 \end{defn}

 \begin{rem}{\label{R:sections}}
  Proposition \ref{P:2catsection} remains valid when $\mfA$,
  $\mfC$ and $\mfD$ are 2-groups, and $D_p(f_0,f)$ is replaced by
  $D_p(f_0,f)_*$. However, the natural functor $D_p(f_0,f)_* \to \mfA$
  is  not in general a fibration anymore.
 \end{rem}

Finally,  observe that the nerve functor $N \: \mathbf{2Cat} \to
\mathbf{SSet}$ respects the difference construction. That is,
$ND_p(f_0,f)$ is naturally homotopy equivalent to
$D_{Np}(Nf_0,Nf)$. This is because $N$ preserves fiber products
(Appendix, Proposition \ref{P:nerve}) and path spaces (Appendix,
Proposition \ref{P:homspace}).

\subsection{Difference fibrations for crossed-modules}{\label{SS:Difference2}}

We saw in the previous  subsection  that we can perform difference
construction for 2-groups. We will make this more explicit using the
language of crossed-modules. So assume $\bbF$, $\bbG$ and $\bbH$ are
crossed-modules, as in the following picture:
    $$\xymatrix@=30pt@M=6pt{         &  \bbG \ar[d]^p  \\
         \bbH  \ar@/^/[ru]|-{f'} \ar@/_/[ru]|-{f}  \ar[r]  &  \bbF    }$$
To avoid notational complications, we have used $f'$ instead of
$f_0$. In what follows, it would be helpful to think of elements of
$G_1$ (respectively, $G_2$) as 1-cells (respectively, 2-cells) of
the nerve $N\bbG$ (see Appendix).

Let $[D_2\to D_1]$ be the crossed-module presentation of
$D_p(f',f)_*$. We will write down exactly what $D_1$ and $D_2$
are. By definition, we have
   $$D_1=\left\{(h,\alpha) \ | \ h \in H_1, \alpha \in \Ker p_2,
             \ \text{s.t.} \ f'_1(h)\underline{\alpha}=f_1(h)\right\}.$$
It should be clear what this means: $\alpha$ is 2-cell that is
vertical (because it is in $\Ker p_2$) and joins the 1-cells
$f'_1(h)$ to $f_1(h)$. The group multiplication is
      $$(h,\alpha)(k,\beta)=(hk,\alpha^{f_1'(k)}\beta).$$

To determine $D_2$, we have to pick a 2-cell $\beta \in H_2$, and
find all pointed vertical homotopies from  $f'_2(\beta)$ to
$f_2(\beta)$. A pointed homotopy from $f'_2(\beta)$ to
$f_2(\beta)$ is a 2-cell $\gamma \in G_2$ such that
$f'_2(\beta)\gamma=f_2(\beta)$.  To ensure it is vertical, we need
to have $\gamma \in \Ker p_2$. This, however,  is automatic, since
$p_2f'_2(\beta)=p_2f_2(\beta)$. The conclusion is that, for any
$\beta \in H_2$, there is a unique vertical homotopy from
$f'_2(\beta)$ to $f_2(\beta)$; it is given by
$\gamma=f'_2(\beta)^{-1}f_2(\beta)$. Therefore,
           $$D_2=H_2.$$
The map  $D_2 \to D_1$ is given by
$$H_2=D_2 \ni\ \beta  \ \mapsto \
   \big(\underline{\beta},\ f'_2(\beta)^{-1}f_2(\beta)\big)\ \in D_1.$$
An element $(h,\alpha) \in D_1$ acts on $\beta \in D_2=H_2$
by sending it to $\beta^h$.

There is a natural map of crossed-modules $q \: D_p(f',f)_*  \to
\bbH$ given by $q_1(h,\alpha)=h$ and $q_2=\id$.

\subsection{The special case}{\label{SS:Special}}
We now consider the special case of the difference  construction
that is relevant to Theorem \ref{T:maps}, and show that it recovers
Definition \ref{D:difference2} (see Proposition \ref{P:compatible}).
Suppose we are given a group $\Gamma$, a 2-group $\bbG$, and a
homomorphism $\chi \: \Gamma \to \pi_1\bbG$. We want to study lifts
of $\chi$ to maps $\Gamma \to \bbG$. Keep in mind that here we are
talking about maps in the {\em homotopy category} of 2-groups. So it
would definitely be false to consider 2-group maps $\Gamma \to
\bbG$. One way to handle the situation is to work with {\em weak
maps}. But this is not very convenient. The better way would be to
pick a cofibrant replacement $\bbH$ for $\Gamma$ (see Example
\ref{E:replacement}), and use the fact that a map $\Gamma \to \bbG$
in the homotopy category of 2-groups can be represented by a map of
2-groups $\bbH \to \bbG$, and that the latter is unique up to
pointed transformation.

Recall (Example \ref{E:replacement}) how we construct a cofibrant
replacement for $\Gamma$: we choose a presentation $F/R=\Gamma$,
where $F$ is a free group, and form the crossed-module $\bbH=[R \to
F]$; the natural map $\bbH \to \Gamma$ is then our cofibrant
replacement. Throughout the paper, we fix such a cofibrant
replacement.

Our problem is to study  the difference construction for the
following situation:
$$
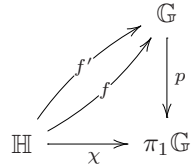
\xymatrix@=30pt@M=6pt{         &  \bbG \ar[d]^p  \\
           \bbH  \ar@/^/[ru]|-{f'} \ar@/_/[ru]|-{f}  \ar[r]_\chi
                                                &  \pi_1\bbG    }$$
All maps are now honest maps of crossed-modules. (We have abused
notation and denoted the induced map $\bbH \to \pi_1\bbG$ also by
$\chi$.) As we saw in the previous subsection, this picture gives a
map of crossed-modules $D_p(f',f)_* \to \bbH$, which is indeed a
fibration of crossed-modules in the sense of Definition
\ref{D:fibrations}. Our aim is to compare this with Definition
\ref{D:difference2}.

Let $(E,\rho)$, respectively $(E',\rho')$, be the object of
$\M(\Gamma,\bbG)$ obtained from pushing out $\bbH$ along $f$,
respectively $f$; see Definition \ref{D:pushout}). Recall
(Definition \ref{D:difference2}) that the difference
$D\big((E',\rho'),(E,\rho)\big)$ is defined to be the following
exact sequence:
   $$1 \to \pi_2\bbG \to E'\times_{\bbG}E \to \Gamma \to 1.$$
We prove the following proposition.

 \begin{prop}{\label{P:compatible}} Notation being as above,
   the map $D_p(f',f)_* \to \bbH$ is a fibration of crossed-modules
   and is naturally equivalent to $E'\times_{\bbG}E \to \Gamma$.
   More precisely, we have the following commutative square in which
   the horizontal arrows are equivalences
   of crossed-modules and the vertical arrows are fibrations:
      $$\xymatrix@=12pt@M=10pt{
            D_p(f',f)_*   \ar[r]^{\sim}\ar[d]  &  E'\times_{\bbG}E \ar[d]  \\
             \bbH \ar[r]^{\sim}        &  \Gamma      }$$
 \end{prop}

 \begin{proof}
    We use the explicit description of the crossed-module
    $D_p(f',f)_*=[D_2\to D_1]$
    given in $\S$\ref{SS:Difference2}:
      $$D_1=\left\{(x,\alpha) \ | \ x \in F, \ \alpha \in G_2, \  \text{s.t.}
         \ f'_1(x)\underline{\alpha}=f_1(x)\right\} \ \ \text{and} \ \ D_2=R.$$
    The map $D_2 \to D_1$ is given by
    $r \mapsto \big(\underline{r},\ f'_2(r)^{-1}f_2(r)\big)$,
    where $\underline{r}$ stands for the image of $r$ in $F$.
    Notice that this is an injection, so we can identify $D_2$ with a subgroup of
    $D_1$.
    Let us now give an explicit description of  $E'\times_{\bbG}E$.
    Recall (Definition \ref{D:pushout}) that
      $$E=F\ltimes^{R,f}G_2   \ \ \text{and} \ \ E'=F\ltimes^{R,f'}G_2.$$
    Using this, we get
       $$E'\times_{\bbG}E=\left\{ \big((x,\alpha),(y,\beta)\big) \ |
          \ \bar{x}=\bar{y}\in \Gamma, \ f'_1(x)\underline{\alpha}=
                               f_1(y)\underline{\beta} \in  G_1\right\}/J,$$
     where $x,y \in F$, $\alpha, \beta \in G_2$,
    and $J$ is the subgroup generated by
          $$\left\{\big((\underline{r},f'_2(r)^{-1}),(1,1\big);
                                        \ \  r \in R \right\}, \ \
           \left\{\big((1,1),(\underline{s},f_2(s)^{-1})\big);
                                        \ \  s \in R \right\},$$
           $$\text{and} \ \ \left\{\big((1,\gamma),(1,\gamma)\big);
                                                    \ \gamma \in G_2 \right\}.$$
    Therefore,
      $$J=\left\{\big((\underline{r},f'_2(r)^{-1}\gamma),
                                  (\underline{s},f_2(s)^{-1}\gamma)\big);
                                    \ \  r,s \in R, \ \gamma \in G_2 \right\}.$$

    Define  the map $\Lambda \: D_1 \to E'\times_{\bbG}E$ by
         $$\Lambda(x,\alpha)=\big((x,\alpha),(x,1)\big).$$
    We claim that $\Lambda$ is surjective and its kernel is $R=D_2$.
    This proves that the map $D_p(f',f)_* \to E'\times_{\bbG}E$
    is an equivalence of crossed-modules.

    \vspace{0.1in}
    \noindent{\em Surjectivity.}  Pick an element
        $a=\big((x,\alpha),(y,\beta)\big)$ in $E'\times_{\bbG}E$.
         Since $\bar{x}=\bar{y}$, there is $s \in R$ such that
        $y=x\underline{s}$. Using the fact that
         $$\big((x,\alpha),(x\underline{s},\beta)\big)
                       = \big((x,\alpha),(x,f_2(s)\beta)\big)$$
        in $E'\times_{\bbG}E$, we may assume that $x=y$, that is,
        $a=\big((x,\alpha),(x,\beta)\big)$.
        On the other hand, after  multiplying on the right by
        $\big((1,\beta^{-1}),(1,\beta^{-1})\big) \in J$,
        we may assume that $\beta=1$; that is
        $a=\big((x,\alpha),(x,1)\big)$.
        This is obviously in the image of $\Lambda$.

    \vspace{0.1in}
    \noindent{\em Kernel of $\Lambda$ is $R$.}        Easy verification.

    \vspace{0.1in}
     To show that
    $D_p(f',f)_* \to \bbH$      is a fibration, we have to show that
    the   map $D_1 \to F$ which sends $(x,\alpha)$ to $x$,
    and the map $\id \: D_2=R \to R$ are surjective (Definition \ref{D:fibrations}).
    The latter is obvious. The former follows from the fact that,
    for every $x \in F$, we have the equality
    $\overline{f_1(x)}=\overline{f_1'(x)}$ in
    $\pi_1\bbG$.
    This is true because both these elements are  equal to $\chi(x)$.

    Commutativity of the square is obvious.
 \end{proof}

\section{Appendix: 2-categories and 2-groupoids}{\label{A:1}}

In this appendix we quickly go over some basic facts and
constructions we need about 2-categories, and fix some terminology.
Most of the material in this appendix can be found in \cite{Notes}.

For us, a {\em 2-category} means a strict 2-category. A {\em
2-groupoid} is a 2-category in which every 1-morphism and every
2-morphism has an inverse (in the strict sense). Every category
(respectively, groupoid) can be thought of as a 2-category
(respectively, 2-groupoid) in which all 2-morphisms are  identity.

A {\em 2-functor} between 2-categories means a strict 2-functor. We
sometimes refer to a 2-functor simply by a functor, or a {\em map of
2-categories}. A 2-functor between 2-groupoids is  simply a
2-functor between the underlying 2-categories.

By {\em fiber product} of 2-categories we mean strict fiber product.
We will not encounter {\em homotopy} fiber product of 2-categories
in this paper.

The terms `morphism', `1-morphism' and `arrow' will be used synonymously.
We use multiplicative notation for elements
of a groupoid, as opposed to the compositional notation (it means, $fg$
instead of $g\circ f$).

\vspace{0.1in} \noindent{\bf Notation.} We use the German letters
$\mfC$, $\mfD$,... for general 2-categories and $\bbG$, $\bbH$,...
for 2-groupoids. The upper case script letters $A$, $B$, $C$... are
used for objects in such 2-categories, lower case script letters
$a$, $b$, $g$, $h$... for 1-morphisms, and lower case Greek  letters
$\alpha$, $\beta$... for 2-morphisms. We denote the category of
2-categories by $\mathbf{2Cat}$ and the category of 2-groupoids by
$\mathbf{2Gpd}$.

\subsection{2-functors, weak 2-transformations
   and modifications}{\label{SA:Transformations}}

We  recall what weak 2-transformations between strict 2-functors
are. More details can be found in \cite{Notes}.  We usually
suppress the adjective weak.

Let $P,Q \: \mathfrak{D} \to \mathfrak{C}$ be (strict) 2-functors.
By a {\em weak 2-transformation} $T \: P \Rightarrow Q$ we mean,
assignment  of an arrow $t_A$ in $\mathfrak{C}$ to every object
$A$ in $\mathfrak{D}$, and a 2-morphism $\theta_c$ in
$\mathfrak{C}$  to every arrow $c$ in $\mathfrak{D}$, as in the
following diagram:
    $$\xymatrix@=12pt@M=10pt{
         P(A) \ar[r]^{t_A}\ar[d]_{P(c)}  &
              Q(A)  \ar[d]^{Q(c)} \ar@{=>}[ld]^{\theta_c}  \\
         P(B)    \ar[r]_{t_{B}}    &  Q(B)   }$$
We require that $\theta_{\id}=\id$, and that $\theta_h$ satisfy the
obvious compatibility conditions with respect to
2-morphisms and composition of morphisms.

A transformation between two weak transformations $T$, $S$,
sometimes called a {\em modification}, is a rule to assign to each
object $A \in  \mfD$
a 2-morphism $\mu_A$ in  $\mfC$ as in the following diagram:
  $$\xymatrix@C=10pt@R=-5pt@M=6pt{ & \ar@{=>}[dd]^{\mu_A}&\\
        P(A) \ar@/^1.1pc/ [rr]^{t_A} \ar@/_1.1pc/ [rr]_{s_A} &&  Q(A)   \\
            & &    }$$
The 2-morphisms $\mu_A$ should satisfy the obvious compatibility
relations with $\theta_c$ and $\sigma_c$, for every arrow $c \: A
\to B$ in $\mathfrak{D}$. (Here $\sigma_c$ are for $S$ what
$\theta_c$ are for $T$.) This relation can be written as
$\sigma_c\mu_A=\theta_c\mu_B$.

 \begin{defn}{\label{D:homo2}}
  Given 2-categories $\mathfrak{C}$ and $\mathfrak{D}$,
  we define the mapping 2-category $\homc(\mathfrak{D},\mathfrak{C})$ to be the
  2-category whose objects are strict 2-functors from
  $\mathfrak{D}$ to $\mathfrak{C}$, whose 1-morphisms
  are weak 2-transformations between 2-functors, and whose
  2-morphisms are modifications. When $\mathfrak{C}$ and $\mathfrak{D}$
  are 2-groupoids, then $\homc(\mathfrak{D},\mathfrak{C})$ is also
  a 2-groupoid.
 \end{defn}

Viewing $\tgp$ as a full subcategory of $\mathbf{2Gpd}$, we can use
the same notion for 2-groups as well. In fact, in the case of
2-groups, we are more interested in the pointed versions of the
above definition. Namely, a {\em pointed} 2-transformation is
required to satisfy the extra condition $t_*=\id$. A pointed
modification is, by definition, the identity modification!

For 2-groups $\bbG$ and $\bbH$, we denote the 2-groupoid of {\em
pointed} weak maps from  $\bbH$ to $\bbG$ by $\homc_*(\bbH,\bbG)$.

\subsection{Nerve of a 2-category}{\label{SA:Nerve}}

We review the nerve construction for 2-categories, and recall its
basic properties \cite{M-S}, \cite{Notes}.

Let $\bbG$ be a 2-category .  We define the {\em nerve} of $\bbG$,
denoted by $N\bbG$, to be the simplicial set defined as follows. The
set of of 0-simplices of $N\bbG$ is the set of objects of $\bbG$.
The 1-simplices are the morphisms in $\bbG$. The 2-simplices are
diagrams of the form
  $$\xymatrix@C=6pt@R=14pt@M=6pt{  & B \ar[dr]^g \ar@{=>}[d]^{\alpha} &   \\
             A \ar[rr]_h \ar[ru]^f   &    &   C   }$$
where $\alpha \: fg \Rightarrow h$ is a 2-morphism. The 3-simplices
of $N\bbG$ are commutative tetrahedra of the form
\label{nervediagram} \label{3cell}
   $$\xymatrix@C=3pt@R=12pt@M=6pt{ &&& D   &&& \\
                               &&&& \ar@{=>}[dl]^{\gamma} && \\
                            &&& B \ar[uu]_l \ar@{=>}[ul]^{\beta}
                                \ar@{=>}[d]^{\alpha} \ar[drrr]^g &  &&   \\
                        A \ar[rrrrrr]_h \ar[urrr]^(0.65)f   \ar[rrruuu]^k
                      && \ar@{:>}[ul]_<{\delta} &&&& C \ar[llluuu]_m   }$$
Commutativity of the above tetrahedron means
$(f\gamma)(\beta)=(\alpha m)(\delta)$. That is, the following
square of transformations is commutative:
            $$\xymatrix@=12pt@M=10pt{
                 fgm  \ar@{=>}[r]^{f\gamma} \ar@{=>}[d]_{\alpha m}
                                          &  fl \ar@{=>}[d]^{\beta}  \\
                            hm     \ar@{=>}[r]_{\delta}        &    k   }$$
For $n\geq 3$, an $n$-simplex of $N\bbG$ is an $n$-simplex such that
each of its sub 3-simplices is a commutative tetrahedron as
described above. In other words, $N\bbG$ is the coskeleton  of the
3-truncated simplicial set $\{N\bbG_0,N\bbG_1,N\bbG_2,N\bbG_3\}$
defined above.

The nerve  gives us a functor $N \: \mathbf{2Cat} \to
\mathbf{SSet}$, where    $\mathbf{SSet}$ is the the category of
simplicial sets.

\subsection{Moerdijk-Svensson closed model structure on 2-groupoids}{\label{SA:M-S}}

We give a quick  review of the Moerdijk-Svensson closed model
structure on the category of 2-groupoids. The main reference is
\cite{M-S}.

 \begin{defn}{\label{D:2fibration}}
     Let $\bbH$ and $\bbG$ be 2-groupoids, and $P \: \bbH \to \bbG$
     a functor between them. We say that
     $P$ is  a {\em fibration}, if it satisfies the
     following properties:
     \begin{itemize}

   \item[$\mathbf{F1.}$] For every arrow $a \: A_0 \to A_1$ in
     $\bbG$, and every object  $B_1$ in $\bbH$ such that $P(B_1)=A_1$,
     there is an object $B_0$ in $\bbH$ and an arrow $b \: B_0 \to B_1$
     such that $P(b)=a$.

   \item[$\mathbf{F2.}$] For every 2-morphism
       $\alpha \: a_0 \Rightarrow a_1$ in
     $\bbG$ and every arrow $b_1$ in $\bbH$ such that $P(b_1)=a_1$,
     there is an arrow $b_0$ in $\bbH$ and a
     2-morphism $\beta \: b_0 \Rightarrow b_1$
     such that $P(\beta)=\alpha$.
  \end{itemize}
 \end{defn}

 \begin{defn}{\label{D:homotopygroup}}
     Let $\bbG$ be a 2-groupoid, and $A$ an object in $\bbG$.
     We define the following.
     \begin{itemize}

   \item $\pi_0\bbG$  is the set of equivalence classes of objects in $\bbG$.

   \item $\pi_1(\bbG,A)$ is the group of 2-isomorphism classes of arrows
         from $A$ to itself. The {\em fundamental groupoid}
         $\Pi_1\bbG$ is  the groupoid whose objects are the same as those
         of $\bbG$ and whose morphisms are 2-isomorphism classes of
         1-morphisms in $\bbG$.

   \item $\pi_2(\bbG,A)$ is the group of 2-automorphisms of the identity
         arrow $1_A \: A \to A$.
  \end{itemize}
  These invariants are functorial with respect to 2-functors.
  A map $\bbH \to \bbG$ is called a {\em (weak) equivalence of 2-groupoids}
  if it induces a bijection on $\pi_0$, $\pi_1$ and $\pi_2$, for every
  choice of a base point.
 \end{defn}

Having defined the notions  of fibration and equivalence between
2-groupoids, we define {\em cofibrations} using the left lifting
property. There is a more explicit description of cofibrations
which can be found in (\cite{M-S} page 194), but we skip it here.

 \begin{thm}[\oldcite{M-S}, Theorem 1.2]{\label{T:QuillenStr}}
  With weak equivalences, fibrations and cofibrations defined as above,
  the category of 2-groupoids has a natural structure of a closed model
  category.
 \end{thm}

The nerve functor is a bridge between the homotopy theory of
2-groupoids and the homotopy theory of simplicial sets. To justify
this statement, we quote the following from \cite{M-S}.

 \begin{prop}[see \oldcite{M-S}, Proposition 2.1]{\label{P:nerve}}\indent\par
    \begin{itemize}

     \item[$\mathbf{i.}$] The functor $N\: \mathbf{2Cat} \to \mathbf{SSet}$
      is faithful,  preserves fiber products,
      and sends transformations between 2-functors to simplicial homotopies.

    \item[$\mathbf{ii.}$] The functor $N$ sends a fibration
      between 2-groupoids (Definition \ref{D:2fibration}) to a
      Kan fibration. Nerve of every
      2-groupoid  is a Kan complex.

    \item[$\mathbf{iii.}$] For every (pointed) 2-groupoid $\bbG$ we have
         $\pi_i(\bbG)\cong\pi_i(N\bbG)$, $i=0,1,2$.

    \item[$\mathbf{iv.}$] A map $f \: \bbH \to \bbG$ of 2-groupoids is an
      equivalence if and only if $Nf \: N\bbH \to N\bbG$ is a weak
      equivalence of simplicial sets.
   \end{itemize}
 \end{prop}

 \begin{rem}{\label{R:nerve2gp}}
  We can think of a 2-group as a 2-groupoid with one object. This
  identifies $\tgp$ with a full subcategory
  of $\mathbf{2Gpd}$. So, we can talk about nerves of 2-groups.
  This is a functor $N \: \tgp \to \mathbf{SSet}_*$,
  where $\mathbf{SSet}_*$ is the category of pointed simplicial sets.
  The above proposition remains valid if we replace 2-groupoids by
  2-groups and $\mathbf{SSet}$ by $\mathbf{SSet}_*$ throughout.
 \end{rem}

The  functor $N\: \mathbf{2Gpd} \to \mathbf{SSet}$ has a left
adjoint $W \: \mathbf{SSet} \to \mathbf{2Gpd}$, called the {\em
Whitehead 2-groupoid} whose definition can be found in (\cite{M-S}
page 190, Example 2).

It is easy to see that $W$ preserves homotopy groups. In particular,
it sends weak equivalences of simplicial sets to equivalences
of 2-groupoids. Much less obvious is the following

 \begin{thm}[\oldcite{M-S},  $\S$2]{\label{T:QuillenEq}}
   The pair
     $$W\:\mathbf{SSet} \rightleftharpoons \mathbf{2Gpd}:N$$
   is a Quillen pair. It satisfies the following properties:
   \begin{itemize}

    \item[$\mathbf{i.}$] Each adjoint preserves weak equivalences.

    \item[$\mathbf{ii.}$] For every 2-groupoid $\bbG$,  the counit
       $WN(\bbG) \to \bbG$ is a weak equivalence

    \item[$\mathbf{iii.}$] For every simplicial set $X$ such that
      $\pi_iX=0$, $i\geq 3$, the unit of adjunction
      $X \to NW(X)$ is a weak equivalence.
   \end{itemize}
   In particular, the functor
   $N \: \Ho(\mathbf{2Gpd}) \to \Ho(\mathbf{SSet})$ induces an
   equivalence of categories between $\Ho(\mathbf{2Gpd})$ and
   the category of homotopy 2-types. (The latter
   is defined to be the full subcategory of $\Ho(\mathbf{SSet})$
   consisting of all $X$ such that $\pi_iX=0$, $i \geq 3$.)
 \end{thm}

 \begin{rem}{\label{R:pointed}}
  The pointed version of the above
  theorem is also valid.  The proof is just a
  minor modification of the proof of the above theorem.
 \end{rem}

The following Proposition follows from Theorem \ref{T:QuillenEq} and
Remark \ref{R:pointed}.

 \begin{prop}{\label{P:equiv}}
   The functor $N \: \Ho(\tgp) \to \Ho(\mathbf{SSet}_*)$
   induces an equivalence between $\Ho(\tgp)$
   and the full subcategory of $\Ho(\mathbf{SSet}_*)$ consisting
   of connected pointed homotopy 2-types.
 \end{prop}

It is also well-known that the geometric realization functor $|-| \:
\mathbf{SSet}_* \to \mathbf{Top}_*$     induces an equivalence of of
homotopy categories. So we have  the following

 \begin{cor}{\label{C:equiv}}
  The functor $|N(-)| \: \Ho(\tgp) \to \Ho(\mathbf{Top}_*)$
  induces an equivalence between $\Ho(\tgp)$
  and the full subcategory of $\Ho(\mathbf{Top}_*)$ consisting
  of connected pointed homotopy 2-types.
 \end{cor}

The following proposition says that derived mapping 2-groupoids
have the correct homotopy type. We have denoted the category
$\{0\to1\}$ by $\mathbf{I}$.

 \begin{prop}{\label{P:homspace}}
  Let $\bbG$ and $\bbH$ be 2-groupoids. Then there is a natural
  homotopy equivalence
       $$N\rhom(\bbH,\bbG)\simeq\Hombf(N\bbH,N\bbG),$$
  where the left hand side is  defined to be $N\homc(\bbF,\bbG)$,
  where $\bbF$ is a cofibrant replacement for $\bbH$, and the right
  hand side is the simplicial  mapping space. In particular,
  $N\rhom(\mathbf{I},\bbG) \simeq (N\bbG)^{\Delta^1}$, that is, the
  nerve of the path category is naturally homotopy equivalent to the
  path space of the nerve.
 \end{prop}

We also have the pointed version of the above proposition.

 \begin{prop}{\label{P:htpyequivalent}}
   Let $\bbG$ and $\bbH$ be 2-groups. Then, there is a natural
   homotopy  equivalence
    $$N\rhom_*(\bbH,\bbG)\simeq\Hombf_*(N\bbH,N\bbG).$$
 \end{prop}

 \begin{rem}{\label{R:weak}}
  In the definition of $\homc$ we have used {\em strict} 2-functors
  as objects, but the arrows are {\em weak} 2-transformation.
  There is a variant of this in which the 2-functors are
  also weak, and this leads to different simplicial mapping spaces
  between 2-group(oid)s. If in  Propositions \ref{P:homspace} and
  \ref{P:htpyequivalent} above we used this simplicial mapping space
  on the left hand side, we would not need to use a cofibrant
  replacement on $\bbH$ to compute the derived mapping spaces.
  In other words, the derived mapping spaces would coincide with
  the actual mapping spaces.
 \end{rem}

\providecommand{\bysame}{\leavevmode\hbox
to3em{\hrulefill}\thinspace}
\providecommand{\MR}{\relax\ifhmode\unskip\space\fi MR }
\providecommand{\MRhref}[2]{%
  \href{http://www.ams.org/mathscinet-getitem?mr=#1}{#2}
} \providecommand{\href}[2]{#2}

\end{document}